\documentclass[10pt,amssymb,amscd]{amsproc}
\usepackage{color}
\usepackage{amssymb}
\usepackage{amsmath}
\usepackage{amsthm} 
\usepackage{mathrsfs}                                                            
\usepackage[all]{xy}
\usepackage{wrapfig}

\def\2{{1\over 2}}

\newcommand{\rf}[1]{(\ref{#1})}

\newcommand{\ud}{\mathrm{d}}
\renewcommand{\t}{\tilde}
\newcommand{\p}{\partial}

\def\p{\partial}

\def\b{\bar}

\def\<{\langle}
\def\>{\rangle}
\def\+{\dagger}

\begin{document}
\title[Representations of the affine $ax+b$-group, $\widehat{sl}(2,\mathbb{R})$]{On the unitary representations of the affine $ax+b$-group, $\widehat{sl}(2,\mathbb{R})$ and their relatives}
\author{Anton M. Zeitlin}
\address{  \newline Department of Mathematics,\newline
Columbia University,
\newline 2990 Broadway, New York,\newline
NY 10027, USA;\newline
IPME RAS, \newline
V.O. Bolshoj pr., 61, 199178, \newline
St. Petersburg \newline
\newline
zeitlin@math.columbia.edu\newline
http://math.columbia.edu/$\sim$zeitlin \newline
http://www.ipme.ru/zam.html  }

\maketitle

\begin{abstract}
This article focuses on two related topics: unitary representations of the loop $ax+b$-group and their relation to a loop version of the $\Gamma$-function 
and the construction of continuous series for the $\widehat{sl(2,\mathbb{R})}$-algebra. Mainly this is a survey of some results obtained by the author and in collaboration with I.B. Frenkel, alongside with the motivation for them from the physical and mathematical points of view.
\end{abstract}

\section{Introduction}
Recently, the continuous series of unitary representations of $U_q(sl(2,{\mathbb R}))$ were constructed \cite{schmudgen}, \cite{faddeev}, 
so that this set of representations is closed under the tensor product \cite{teschner}, generating a "continuous" tensor category. 
These representations are constructed via the representations of  quantum plane and do not have the classical limit. The representation theory of the quantum plane is very subtle (see e.g. \cite{ivan}, \cite{hyunkyu}) and is related to the properties of quantum dilogarithm. In \cite{ivan} it was shown that the representations of quantum plane and 
furthermore, the intertwining operators of the tensor product representation have a classical limit, where they correspond to appropriate objects in the representation theory of the $ax+b$-group, the affine group of the line (below it is denoted as $G$). In the work of I.B. Frenkel and H.-K. Kim
\cite{hyunkyu} it was shown that basic structures from the quantum Teichmueller theory \cite{chekhov} can be derived from the tensor products  of the representations of the quantum plane, which is very important for the proper understanding of the Chern-Simons theory associated with $SL(2,{\mathbb R})$ (see below for more details).

One important problem to think about is the construction of the affine analogue of the above formalism. Namely, one of the main problems for investigation would be to generalize the Kazhdan-Lusztig equivalence of categories \cite{KL}, \cite{efk} to the case of the quantum plane and $U_q(sl(2,{\mathbb R}))$. It was conjectured in \cite{fzsl2} by I.B. Frenkel and the author, that 
there exists a braided tensor category of representations of the 
affine algebra $\widehat{sl}_k(2,\mathbb{R})$, equivalent to the braided tensor category of the mentioned above continuous series of $U_q(sl(2,\mathbb{R}))$.

One fact supporting this conjecture is that there exists a category of representations 
of the Virasoro algebra related to Liouville theory \cite{teschner2},\cite{teschner3}, which is equivalent to the mentioned category of representations of the Virasoro algebra. Therefore, the corresponding category of representations of the affine algebra $\widehat{sl}_k(2,\mathbb{R})$ is a missing piece in this equivalence of categories.

However, a first natural problem to solve is to find unitary representations for the $\widehat{ax+b}_k$-group, the loop version of the $ax+b$-group with central extension (with central charge $k$) for the $a$-subgroup, and compare them and their tensor structure with the representation theory of quantum plane.

In \cite{zeit} the representation theory of the loop counterpart of the 
$ax+b$-group ($\Omega G$) and its central extension ($\hat G$) were studied. There the unitary representations of $\hat{G}$ are constructed, which naturally generalize the irreducible representations of $G$. 

It is known that the irreducible unitary representations of $G$ divide into three classes (up to equivalence). The first two classes contain infinite-dimensional representations and there is only one representation in each class, however the third class consists of 1-dimensional representations, labeled by the real parameter.
 The infinite dimensional representations of $G$ can be realized in the Hilbert space $L^2({\mathbb R}_{+},\frac{dx}{x})$.  These representations are produced by the canonical construction 
of induced representations. This classical material is surveyed in Section 2.1.

However, in the case of the loop $ax+b$-group the 
situation is more subtle. In \cite{zeit} the $L^2$ space with respect to the 
Wiener measure is considered as a representation space of $\hat{G}$. 
The constructed unitary representations of $\hat{G}$ are labeled 
by a certain function. It is proven in \cite{zeit} that certain representations in this class (e.g. when this function is  
constant) are irreducible. We review these results in Section 3.

An interesting issue studied in \cite{zeit} is based on the mentioned relations of 
the representations of group $G$ and the theory of special functions. It is known that the representations of $G$ are related to the $\Gamma$- and $B$- functions (see e.g. \cite{vilenkin}). Therefore, one may hope to construct their loop counterparts while studying the associated loop group. In \cite{zeit} the first step in this direction is made. It is known, that one can relate (via the bilateral Laplace transform) the 
action of the group element of representation of $G$ in the representation space with the integral operator, so that its  kernel is expressed via the $\Gamma$-function. 
In \cite{zeit} the loop analogue of the $\Gamma$-function was defined, considering similar arguments applied to $\hat G$. It is a functional on the space of paths,   
which we denote as 
$\hat{\Gamma}_{\mu}$, depending on some continuous function $\mu$, such that $\mu(u)>0$ on $[0,2\pi]$, and has a property which generalizes the famous property of the $\Gamma$-function, $\Gamma(x+1)=x\Gamma(x)$:
\begin{eqnarray}
&&\int_0^{2\pi} g(v)\mu(v)\hat\Gamma_{\mu}(z+\delta_v)dv=\int _{0}^{2\pi}g(v)z(v)dv\hat\Gamma_{\mu}(z)+\nonumber\\
&&\frac{1}{t}\int_0^{2\pi}g''(v)\frac{\delta}{\delta z(v)}\hat\Gamma_{\mu}(z)dv, 
\end{eqnarray}
where $g(v)$ is any twice differentiable function on $[0,2\pi]$, such that $g(0)=g(2\pi)=0$,  
$\delta_v=\delta(u-v)$ is a delta-function on the interval $[0,2\pi]$ and $t$ is the parameter of the Wiener measure. 
We describe this construction of the loop $\Gamma$-function from \cite{zeit} in Section 4. 

It is known (see e.g. \cite{ivan}) that the unitary representations of $G$ are closed under the tensor product and there are three types of ''simple'' objects in the category of unitary representations of $G$. It appears that their tensor products decompose as direct integrals of these ''simple'' objects. 
Using the principles discussed in the beginning of this section, it is expected to have the braided tensor category for 
$\hat{G}$, where the braiding is related to the value of the central charge. One can also hope to obtain a differential equation governing the intertwining operators, i.e. the analogue of the Knizhnik-Zamolodchikov equation. 

However, the results related to the $\widehat{ax+b}_k$-group are only part of the main task, namely the construction of the continuous series of unitary representations of $\widehat{sl}_k(2,{\mathbb R})$. 

Unfortunately, the standard approach of inducing representation of $\widehat{sl}_k(2,{\mathbb R})$ wouldn't fit the construction, since the resulting modules appear to be nonunitary. 
It turns out one can use the results obtained in \cite{zeit} to construct new modules for $\widehat{sl}_k(2,{\mathbb R})$ by means of the Wakimoto-type formalism, using the "currents" corresponding 
to the Lie algebra elements of $(ax+b)_k$ and infinite dimensional Heisenberg algebra  free fields. This was done in \cite{fzsl2}, inspired by the formalism of \cite{zhu}. It turned out that the correlators of the resulting $\widehat{sl}_k(2,{\mathbb R})$-currents, defining the pairing in the  representation,  diverge, and therefore we had to describe the scheme of eliminating those divergencies. This leads to a very interesting graphical formalism, similar to the Feynman diagram technique, where the divergences corresponded to 1-loop graphs. The regularization scheme involves dependence on the infinite family of parameters: one parameter for each loop with a given number of vertices.
These modules were called {\it continuous series} for $\widehat{sl}_k(2,{\mathbb R})$ in \cite{fzsl2}, since the described above Wakimoto-type realization is a natural affinization of certain realization of standard continuous series of $sl(2,\mathbb{R})$. These representations are labeled by the continuous parameter arizing from the highest weight of the Fock module of the infinite-dimensional Heisenberg algebra. 
One can generalize this construction of continuous series to the higher rank case, using similar procedure, that works for quantum groups, introduced in \cite{fip}. 
We describe the construction of the above $\widehat{sl}_k(2,{\mathbb R})$ representations from \cite{fzsl2} in Sections 5 and 6 of this paper. Section 5 is devoted to the construction of Wakimoto-like modules, inspired by the construction of the standard continuous series of $sl(2,\mathbb{R})$, sketched in Sections 2.2 and 2.3. Section 6 is devoted directly to the  construction of the $\widehat{sl}_k(2,{\mathbb R})$-modules of \cite{fzsl2} using regularization and graphical calculus.

An important question one can immediately ask is for which values of regularization parameters the resulting representations are unitary. 
Answering this question will involve studying the resulting bilinear form and analyzing it using the graph formalism that was a cornerstone of the definition of those representations in \cite{fzsl2}. 
However, we do believe that finding the representations forming a braided tensor category will single out the necessary family  of unitary representations. Therefore the next question to ask is what are the intertwining operators for the proposed tensor category. One of the ways to do that is to construct them in a similar fashion as in the Virasoro case \cite{teschner2}, \cite{teschner3}.  

There is an interesting physical perspective on the proposed tensor categories discussed, since currently the study of Chern-Simons theories with noncompact gauge groups has become very important from the point of view of both mathematics and physics, see e.g. \cite{dimofte} for a review. In particular, it was argued since the late 80s  \cite{verlinde} that the canonical quantization of $SL(2,{\mathbb R})$ Chern-Simons theory is connected with both the quantum Teichmueller theory \cite{chekhov}, which is related to the representation theory of the quantum plane and the Liouville theory, related to representation theory of $U_q(sl(2,{\mathbb R}))$. It is known that the Liouville theory can be obtained from the $SL(2,{\mathbb R})$ WZW theory by means of the Drinfeld-Sokolov reduction.  At the same time, the study of Chern-Simons theory with the compact gauge group $G$ showed, that its space of states in the presence of Wilson lines is isomorphic to the space of conformal blocks of the WZW model associated with $G$ \cite{witten}. Therefore, we may expect that this works for noncompact Lie groups too, and the construction of the continuous series of unitary representations of $\widehat{sl}_k(2,{\mathbb R})$ will provide an important link between the quantum Teichmueller theory, the $SL(2,{\mathbb R})$ Chern-Simons, the $SL(2,{\mathbb R})$ WZW model and its reduction, the Liouville theory.

In the end of this introduction let us summarize shortly for convenience the contents of this paper. In Section 2 we gather some standard facts about the unitary representations of $ax+b$-group and  $sl(2,\mathbb{R})$ Lie algebra. In Section 3, we describe some  unitary representations of the loop $ax+b$-group, which are generalizations of the ones considered in Section 2. Section 4 is devoted to the construction of a loop analogue $\Gamma$-function, which appears naturally from the representations of loop $ax+b$-group. Section 5 is devoted to the study of auxiliary current algebras, which we will use in Section 6 to construct representations of 
$\widehat{sl}(2,\mathbb{R})$.\\

\noindent{\bf Acknowledgements.} 
I am indebted to Igor Frenkel for proposing the areas described in this review and for our fruitful and pleasant collaboration. I am very grateful to the organizers of the Southeastern Lie Theory Workshop for inviting me with a talk and stimulating discussions. I am indebted to the referee for her/his comments and to Antonina N. Fedorova for careful reading of the manuscript. 

\section{Unitary representations of the affine group of a line and the  related realizations of $sl(2,\mathbb{R})$}

\noindent{\bf 2.1. $ax+b$-group and $\Gamma$-function.} In this subsection we remind the reader of all the necessary facts (for more information and references one can consult \cite{vilenkin}, chapter 5) about the unitary representations of $ax+b$-group, which is the affine group of the real line. In other words, each group element $g=g(a,b)$ is determined by the pair of real numbers $a,b$ such that $a>0$. The composition law is defined as follows: $g(a_1,b_1)g(a_2,b_2)=g(a_1a_2, a_1b_2+b_1)$. In the following we will call this group $G$.

The unitary representations of this group are constructed by means 
of the method of induced representations. Let us consider the representation $r_{\lambda}$ of the subgroup $B$ of $G$ generated by $g(1,b)$-elements, such that $r_{\lambda}(g(1,b))=e^{\lambda b}$ (here $\lambda$ is complex). Then according to the construction of induced representations we have to consider the space of complex valued functions on $G$, i.e. $f(g(a,b))\equiv f(a,b)$ such that
\begin{eqnarray}
f(a,b+b_0)=e^{\lambda b_0} f(a,b).
\end{eqnarray}
Therefore $f(a,b)= e^{\lambda b}f(a,0)$ i.e. the function $f$ can be expressed in terms of the function $\phi(a)=f(a,0)$ on the subgroup $A$ generated by $g(a,0)$-elements. 
Then the operator of induced representation $R_{\lambda}$ acts on $f(g)$ as $R_{\lambda}(g_0)f(g)=f(gg_0)$, and the resulting formula for representation on the functions $\phi(a)$ is:
\begin{eqnarray}
R_{\lambda}(g_0)\phi(a)=e^{\lambda a b_0}\phi(a_0a),
\end{eqnarray}
or, in other words, since the space of functions $\phi$ is just a space of functions on the ray $0<x<\infty$, we have
$R_{\lambda}(g)\phi(x)=e^{\lambda bx}\phi(ax)$. One can consider the invariant measure on $A$, it is 
just $\frac{dx}{x}$. Then the following theorem is valid \cite{vilenkin}.\\

\noindent{\bf Theorem 2.1.} {\it i) Representation $R_{\lambda}$ of $G$ is unitary on $L^2(\mathbb{R}_{+},\frac{dx}{x})$ if $\lambda\in i\mathbb{R}$.\\
ii) Representations $R_{\lambda}$ and $R_{e^{\xi}\lambda}$ are equivalent if $\xi\in \mathbb{R}$. \\
iii) Consider a semigroup $G_+$(resp. $G_-$), consisting of such $g(a,b)$, that $b>0$ (resp. $b<0$). Then $R_{\lambda}$ with $\lambda\in \mathbb{C}$, such that $Re\lambda<0$ (resp. $Re\lambda>0$) is a representation for $G_+$(resp. $G_-$) on $L^2(\mathbb{R}_{+},\frac{dx}{x})$.}\\

As a consequence, one can see that there are only three classes of unequivalent unitary representations $R_{\lambda}$: $R_{\pm i}$ and $R_0$. The representations $R_{\pm i}$ are irreducible, while $R_0$ decomposes into the direct integral of 1-dimensional representations $T_{\rho}$, such that $T_{\rho}(g(a,b))=a^{i\rho}$. One can show that any other irreducible representation of $G$ falls into one of the classes $R_{\pm i}$,  $T_{\rho}$.

In order to generalize representations $R_{\lambda}$ to the loop case, it is useful to consider another form of the representations $R_{\lambda}$. Namely, instead of the space $L^2(\mathbb{R}_{+},\frac{dx}{x})$ one can 
consider $L^2(\mathbb{R},dt)\equiv L^2(\mathbb{R})$, by substitution $x=e^t$. Therefore, the formula for the representation $R_{\lambda}$ can be rewritten as follows:
\begin{eqnarray}\label{rl}
R_{\lambda}(g(e^{\alpha},b))f(t)=e^{i\lambda be^t}f(t+\alpha),
\end{eqnarray}
where $f \in L^2(\mathbb{R},dt)$ and we represented $a$ as $e^{\alpha}$. 

Now we discuss the relation of the representation $R_{\lambda}$ and $\Gamma$-function. It is well known, that the  
bilateral Laplace transform (or the Fourier transform in the complex domain)
\begin{eqnarray}
\mathcal{L}f(p)=\frac{1}{\sqrt{2\pi}}
\int_{\mathbb{R}}e^{ipt}f(t)dt,
\end{eqnarray}
where $p$ is a complex number, has the inverse:
\begin{eqnarray}
\mathcal{L}^{-1}g(t)=\frac{1}{\sqrt{2\pi}}
\int_{\mathbb{R}+iT}e^{-ipt}g(p)dp,
\end{eqnarray}
where $T$ is a real number, so that the contour of integration is in the region of convergence of $g(p)$.
Therefore, if one can make sense of  $\mathcal{L}R_{\lambda}\mathcal{L}^{-1}$, it gives us a representation, equivalent 
to $R_{\lambda}$. Let $\mathcal{D}\subset L^2(\mathbb{R})$ be the space of $C^{\infty}$ functions with finite support. \\

\noindent{\bf Proposition 2.1.} \cite{vilenkin} {\it The bilateral Laplace transform of an element of $\mathcal{D}$ is an analytic function in the entire complex plane of the exponential type, i.e. $|\mathcal{L}f(x+iy)|<Ce^{b|y|}$, where $C>0$, 
$b>0$. At the same time, the inverse Laplace transform of such an  analytic function belongs to $\mathcal{D}$. Moreover, we have the following property:
\begin{eqnarray}
\int_{\mathbb{R}}|\mathcal{L}f(x+iy)|^2dx<\infty,
\end{eqnarray}
where $f\in \mathcal{D}$.}\\

\noindent We notice that $\mathcal{D}$ is invariant under the action of the operators $R_{\lambda}(g)$. 
Let $\Gamma(z)$ denote the $\Gamma$-function: 
\begin{eqnarray}
\Gamma(z)=\int^{\infty}_0 e^{-x}x^{z-1}dx.
\end{eqnarray}
Then we have a Proposition.\\

\noindent {\bf Proposition 2.2.} \cite{vilenkin} {\it i) Consider the action of the $\mathcal{L}R_{\lambda}(g)\mathcal{L}^{-1}$ on $\mathcal{L}\mathcal{D}$, when $g=g(a,b)\in G_+$ and $\lambda<0$. Then 
\begin{eqnarray} 
\mathcal{L}R_{\lambda}(g)\mathcal{L}^{-1}f(t_1)=
\frac{1}{2\pi } \int_{\mathbb{R}+i0}\Gamma(it_1-it_2)a^{-it_1}\Big(\frac{-\lambda b}{a}\Big)^{it_2-it_1}f(t_2)dt_2,
\end{eqnarray}
where $f$ is an analytic function on $\mathbb{C}$ of the exponential type. \\
ii) One can analytically continue the expression above with respect to $-\lambda b$ to all complex plane except for the 
negative real axis. }\\

\noindent {\bf 2.2. $\mathcal{A}$ and $\mathcal{K}$ algebras and their unitary representations.} 
In the rest of this section we present a construction of continuous series of $sl(2,\mathbb{R})$ in a special way, which will be convenient for generalizations to the loop case. First, we construct certain modifications of the Lie algebra of $ax+b$-group, which we call 
$\mathcal{A}$-algebra and $\mathcal{K}$-algebra: they have a similar algebraic structure, but different *-structure.
The $\mathcal{A}$-algebra is an algebra with three generators: $h,e^{\pm}$, so that  $ih,ie^{+}$ generate the Lie algebra of 
ax+b-group, however we make the b-subgroup generator $e^{+}$ to be invertible, so that the inverse element is $e^{-}$.
The commutation relations and the star structure are as follows:
\begin{eqnarray}
[h, e^{\pm}]=\pm ie^{\pm}, \quad e^{\pm}e^{\mp}=1, \quad h^*=h, \quad e^{{\pm}^*}=e^{\pm}.
\end{eqnarray}
The generators $h, \alpha^{\pm}$ of $\mathcal{K}$-algebra have similar commutation relations, but the star-structure is different:
 \begin{eqnarray}
 [h, \alpha^{\pm}]=\mp \alpha^{\pm}, \quad \alpha^{\pm}\alpha^{\mp}=1, \quad h^*=h, \quad {\alpha^{\pm}}^*=\alpha^{\mp}.
\end{eqnarray}
To construct the unitary representations of the above *-algebras, one can realize the generators $h,e^{\pm}$ as unbounded self-adjoint operators $i\frac{d}{dx}, e^{\pm x}$ in the Hilbert space of square-integrable functions on the real line $L^2(\mathbb{R})$ as well as generators $h, \alpha^{\pm}$ as the operators $i\frac{d}{d\phi}, e^{\pm i\phi}$ in the Hilbert space of square integrable functions on the circle $L^2(S^1)$. These generators acting on a specific vector 
in the appropriate Hilbert space generate a dense set. More explicitly one has the following proposition (see e.g. \cite{fzsl2}, Proposition 2.1).\\ 

\noindent {\bf Proposition 2.3.} {\it i) Let $D_A$ be the space spanned by the vectors $a\cdot v_0$, where $v_0=e^{-t x^2}\in L^2(\mathbb{R})$, $t>0$ and 
$a$ belongs to the universal enveloping algebra of $\mathcal{A}$-algebra, such that the action of generators is realized as   
$h=i\frac{d}{dx}, e^{\pm}= e^{\pm x}$. Then $D_A$ is a dense set in $L^2(\mathbb{R})$.\\

\noindent ii) Let $D_K$ be the space spanned by the vectors $a\cdot v_0$, where $v_0=1\in L^2(S^1)$ and 
$a$ belongs to the universal enveloping algebra of $\mathcal{K}$-algebra, such that the action of generators is realized as   
$h=i\frac{d}{d\phi}, \alpha^{\pm}= e^{\pm i \phi}$, where $\phi\in [0,2\pi]$ is the angle coordinate on $S^{1}$. Then $D_K$ is a dense set in $L^2(S^1, d\phi)\equiv L^2(S^{1}).$}\\

\noindent {\bf 2.3. Continuous series of $sl(2,\mathbb{R})$ via $\mathcal{A}$ and $\mathcal{K}$ algebras.} 
 First of all, let us introduce some notations. We are interested in unitary representations of $sl(2,\mathbb{R})$ algebra, i.e. Lie algebra with generators $E,F, H$ such that 
 \begin{eqnarray}\label{sl2r}
&& [E,F]=H, \quad [H, E]=2E,\quad [H, F]=-2F\nonumber\\
&& E=-E^*, \quad F=-F^*, \quad H=-H^* .
 \end{eqnarray}
We remind that there are two standard realizations of continuous series of $sl(2,\mathbb{R})$, one is related to inducing the representations from the diagonal subgroup of $sl(2,\mathbb{R})$, corresponding to the generator $H$, the other one is related to  inducing representations from maximally compact subgroup generated by $J^3=E-F$. It is convenient to introduce the following change of generators (which provides the correspondence between $sl(2,\mathbb{R})$ and $su(1,1)$):
\begin{eqnarray}
J^{\pm}=E+F\mp i H ,  
\end{eqnarray}
so that the relations \rf{sl2r} can be rewritten in the following way:
\begin{eqnarray}
&&[J^3, J^{\pm}]=\pm 2iJ^{\pm}, \quad 
[J^+, J^-]=-iJ^3, \nonumber \\
&&{J^+}^*=-J^-, \quad {J^3}^*=-J^3.
 \end{eqnarray}
Now we will give two realizations of $sl(2,\mathbb{R})$ via generators of $\mathcal{A}$- and $\mathcal{K}$-algebras.
Using $\mathcal{A}$-algebra, one can write down the expressions for $E, F, H$:
\begin{eqnarray}\label{real1}
&&E=\frac{i}{2}(e^+h+he^+)+i\lambda e^+,\nonumber\\
&&F=-\frac{i}{2}(e^-h +he^-)+i\lambda e^-,\\
&&H=-2ih,\nonumber
 \end{eqnarray}
where $\lambda$ is a real parameter. 
Similarly, using $\mathcal{K}$-algebra, one can represent $J^3, J^{\pm}$ in a similar way:
\begin{eqnarray}\label{real2}
&&J^+=\frac{i}{2}(\alpha^+h+h\alpha^+)-\lambda \alpha^+\nonumber\\
&&J^-=\frac{i}{2}(\alpha^-h +h\alpha^-)+\lambda \alpha^-\\
&&J^3=2ih\nonumber
 \end{eqnarray}
\noindent Therefore the following theorem is valid (see e.g. \cite{fzsl2}, Theorem 2.1).\\

\noindent {\bf Theorem 2.2.} {\it Let $\mathcal{D}_A$ be the space spanned by the vectors $a\cdot v_0$, where $v_0=e^{-t x^2}\in L^2(\mathbb{R})$ and 
$a$ belongs to the universal enveloping algebra of $sl(2,\mathbb{R})$-algebra, such that the action of generators is realized as in \rf{real1}, so that  
$h=i\frac{d}{dx}, e^{\pm}= e^{\pm x}$. Then $\mathcal{D}_A$ is a dense set in $L^2(\mathbb{R})$ and it is a representation space for $sl(2,\mathbb{R})$.

\noindent ii) Let $\mathcal{D}_K$ be the space spanned by the vectors $a\cdot v_0$, where $v_0=1\in L^2(S^1)$ and 
$a$ belongs to the universal enveloping algebra of $sl(2,\mathbb{R})$-algebra, such that the action of generators is realized as in \rf{real2}, so that $h=i\frac{d}{d\phi}, e^{\pm}= e^{\pm i \phi}$, where $\phi\in [0,2\pi]\cong S^1$. Then $\mathcal{D}_K$ is a dense set in $L^2(S^1)$ and it is representation space for $sl(2,\mathbb{R})$.} \\

\section{Unitary representations of the path and loop versions of $ax+b$- group}

\noindent {\bf Notation.} In the following we will use several functional spaces, so let us fix notations. Let $\mu$ be 
the $\sigma$-additive measure on some space $M$. Then we denote as $L^2(M,d\mu;k)$ the Hilbert space 
of square-integrable complex-valued functions (here $k$ stands for $\mathbb{C}$ or $\mathbb{R}$). Since in the most of cases we will deal with $L^2(M,d\mu;\mathbb{C})$, we drop $\mathbb{C}$, i.e. $L^2(M,d\mu)\equiv L^2(M,d\mu;\mathbb{C})$. 
The space of real-valued continuous functions on the interval $[a,b]$ will be referred to as $C[a,b]$. The space of real-valued absolutely continuous functions on $[a,b]$, i.e. differentiable functions from $C[a,b]$, whose derivative belong to $L^2([a,b],dx;\mathbb{R})$, where $dx$ is the Lebesgue measure on $[a,b]$ will be denoted as $C'[a,b]$. 

We introduce $C_0'[0,2\pi]$, the Hilbert space of real-valued absolutely continuous functions, such that $x(0)=0$ for any $x\in C_0'[0,2\pi]$. The completion of $C_0'[0,2\pi]$ is the Banach space $C_0[0,2\pi]$ of real-valued continuous functions such that $x(0)=0$ for any $f\in C_0[0,2\pi]$ (see Appendix). 
Finally,  let $C_{0,X}[0,2\pi]$ denote the closed subspace of  $C_0[0,2\pi]$, consisting of real-valued continuous functions such that $x(2\pi)=X$ for some $X\in \mathbb{R}$. 
\\

\noindent{\bf 3.1. Path groups and loop groups associated to $G$.} In this subsection, we define all the path and loop groups, which we study in this article. All the necessary facts concerning the Wiener measure can be found in Appendix. 

At first we define the path group $PG$. Let us consider a set of elements of the form $g(e^{\alpha}, b)$, where $\alpha$ and $b$ are real valued absolutely continuous functions on $[0, 2\pi]$. We can define the bilinear operation as follows:
\begin{eqnarray}\label{mult}
g(e^{\alpha_1}, b_1)\cdot g(e^{\alpha_2}, b_2)=g(e^{\alpha_1+\alpha_2}, e^{\alpha_1}b_2+b_1),
\end{eqnarray}
where $e^{\alpha_1}b_2$ stands for pointwise multiplication of the absolutely continuous functions $e^{\alpha_1}$ and $b_2$. It is clear that the operation is well defined and satisfies the group laws. 
Let us denote the group which is a set of all such elements g($e^{\alpha}, b$) as $PG$. 

In the following we will sometimes need $\alpha$ to be in the Cameron-Martin space for the Wiener measure, so it is useful to introduce the based path group $P_0G\subset PG$, generated by $g(e^{\alpha},b)$, where $\alpha(0)=0$.

Similarly, one can define loop groups $\Omega G$ and $\Omega_0 G$. The loop group $\Omega G\subset PG$ is generated by $(e^{\alpha},b)$, where $\alpha, b$ are such that $\alpha(0)=\alpha(2\pi)$, $b(0)=b(2\pi)$, while 
the based loop group $\Omega_0 G\equiv P_0G\cap\Omega G$.

Finally, one can define the central extended versions of $\Omega G$ (resp. $\Omega_0 G$). Let us consider the elements $g(e^{\alpha}, b, s)$, where $s\in \mathbb{R}$ and $(e^{\alpha}, b)\in \Omega G$ (resp. $\Omega_0 G$). 
Then the multiplication law can be modified in the following way:
\begin{eqnarray}
&&g(e^{\alpha_1}, b_1, s)\cdot g(e^{\alpha_2}, b_2, t)=\nonumber\\
&&g(e^{\alpha_1+\alpha_2}, e^{\alpha_1}b_2+b_1, t+s+
k\int_0^{2\pi}\alpha_1(u)\alpha'_2(u)\ud u),
\end{eqnarray}
where $k\in \mathbb{R}$. We denote the corresponding group $\hat{G}$ (resp. $\hat{G}_{0}$) and we will refer to $k$ as central charge.\\

\noindent{\bf 3.2. Unitary representations.} 
In this subsection, we will describe the unitary representations 
of both $PG$ and $\Omega G$, as well as their subgroups $P_0G$ and $\Omega_0 G$, which are the appropriate generalizations of the unitary representations of $G$, considered in Section 2. 

We will start from $P_0G$ and use the same approach as before, i.e. we will use the method of induced representations. Let us take one-dimensional representations of the $B$-subgroup (i.e. subgroup of elements of the 
form $g(1,b)$) of the form $\tilde{r}_{\lambda}(g(1,b))=e^{\int_0^{2\pi}\lambda(u)b(u)\ud u}$, where $\lambda\in L^2[0,2\pi]$. Following the method of induced representations, like we did in the case of group $G$, we arrive to the following formula for the representations on 
the space of functionals on the $A$-subgroup:
\begin{eqnarray}
\tilde{R}_{\lambda}(g(e^{{\alpha}_0},{b_0}))f(\alpha)=e^{\int_0^{2\pi}\lambda(u)b_0(u)e^{{\alpha}(u)}\ud u}f(\alpha+\alpha_0).
\end{eqnarray}
However, in order to make these representations unitary one needs to define the proper inner product on the space 
of functionals. To do that, one has to consider the Wiener measure (see Appendix). It is not invariant under translations, so one should improve the formula for the representations in order to make them unitary.  
Let us consider the Hilbert space $H^0_p=L^2(C_0[0,2\pi], dw^t;\mathbb{C})$. The following statement is true.\\

\noindent {\bf Theorem 3.1.} {\it The following action of $P_0G$ on the space of functionals of continuous functions  
\begin{eqnarray}\label{act}
&&\rho_{\lambda}(g(e^{\alpha},b))(f)(x)=\nonumber\\
&&e^{-\frac{1}{4t}\int_0^{2\pi}\alpha'(u)\alpha'(u)\ud u-\frac{1}{2t} \int_0^{2\pi}\alpha'(u)\ud x(u)}
e^{\int_0^{2\pi} \lambda(u)b(u)e^{x(u)}\ud u}f(x+\alpha)
\end{eqnarray}
defines the unitary representation of $P_0G$ on $H^0_p$ iff $i\lambda(u)\in L^2([0,2\pi];\mathbb{R})$.}\\

\noindent {\bf Proof.} First we prove that this action is actually an action of a group, i.e. 
$\rho_{\lambda}(g(e^{\alpha_1},b_1))\rho_{\lambda}(g(e^{\alpha_2},b_2))=\rho_{\lambda}(g(e^{\alpha_1+\alpha_2}, e^{\alpha_1}b_2+b_1))$. We 
check it, writing the explicit expressions:
\begin{eqnarray}
&&\rho_{\lambda}(g(e^{\alpha_1},b_1))\rho_{\lambda}(g(e^{\alpha_2},b_2))f(x)=\nonumber\\
&&\rho_{\lambda}(g(e^{\alpha_1},b_1))(e^{-\frac{1}{4t}\int_0^{2\pi}\alpha_2'(u)\alpha_2'(u)\ud u-\frac{1}{2t} \int_0^{2\pi}\alpha_2'(u)\ud x(u)}\nonumber\\
&&e^{\int_0^{2\pi} \lambda(u)b_2(u)e^{x(u)}\ud u}f(x+\alpha_2))=\nonumber\\
&&e^{-\frac{1}{4t}\int_0^{2\pi}\alpha_2'(u)\alpha_2'(u)\ud u-\frac{1}{4t}\int_0^{2\pi}\alpha_1'(u)\alpha_1'(u)\ud u-\frac{1}{2t}\int_0^{2\pi}\alpha_1'(u)\alpha_1'(u)\ud u}\nonumber\\
&&e^{-\frac{1}{2t} \int_0^{2\pi}(\alpha_1'(u)+\alpha_2'(u))\ud x(u)}e^{\int_0^{2\pi} \lambda(u)(b_2(u)e^{\alpha_1(u)}+b_1(u))e^{x(u)}\ud u}f(x+\alpha_1+\alpha_2))=\nonumber\\
&&\rho_{\lambda}(e^{\alpha_1+\alpha_2}, e^{\alpha_1}b_2+b_1)f(x).
\end{eqnarray}
Now we prove that the action $\rho_{\lambda}$  defines a unitary representation for $P_0G$. In order to do that we just need to use the definition of unitarity and the translation invariance property of the Wiener measure:
\begin{eqnarray}
&& (\rho_{\lambda}(g(e^{\alpha},b))(f),\rho_{\lambda}(g(e^{\alpha},b)) (h))_{H_p}=\nonumber\\
&&\int_{C_0[0,2\pi]} \overline{\rho_{\lambda}(g(e^{\alpha},b))(f)(x)}\rho_{\lambda}(g(e^{\alpha},b)) (h)(x)dw^t(x)=\nonumber\\
&&\int_{C_0[0,2\pi]}  e^{-\frac{1}{2t}\int_0^{2\pi}\alpha'(u)\alpha'(u)\ud u-\frac{1}{t}\int_0^{2\pi}\alpha'(u)\ud x(u)}\overline{f(x+\alpha)}g(x+\alpha)dw^t(x)=\nonumber\\
&&\int \overline{f(x)}g(x)dw(x)=(f,g)_{H_p}.
\end{eqnarray}
Thus the theorem is proven.\hfill$\blacksquare$\\

\noindent In the following we will denote representations $\rho_{\lambda}(g(e^{\alpha},b))\equiv \rho_{\lambda}(e^{\alpha},b)$. 
 
 Using the same steps as in the proof of the theorem above, one can show that analogous fact for the Hilbert space $H^{0,X}_p=L^2(C_{0,X}[0,2\pi], dw_X^t;\mathbb{C})$ and group $\Omega_0 G$ is true.\\

\noindent {\bf Proposition 3.1.} {\it i) Formula \rf{act} defines a unitary representation $\rho^X_{\lambda}$ of $\Omega_0 G$ on $H^{0,X}_p$.\\
ii) Representation $\rho^{X_1}_{\lambda}$ on $H^{0,X_1}_p$ is equivalent to the representation of $\rho^{X_2}_{e^{\eta}\lambda}$ on $H^{0,X_2}_p$, where $\eta$ is absolutely continuous, $\eta(0)=0$, $\eta(2\pi)=X_1-X_2$.}\\

\noindent Therefore, all representations $\rho^X_{\lambda}$ are equivalent to the representations $\rho^0_{e^{\xi}\lambda}$ on 
$H^0_l\equiv L^2(C_{0,0}[0,2\pi], dw^t;\mathbb{C})$, where $\xi$ is absolutely continuous, $\xi(0)=0, \xi(2\pi)=X$. 
Moreover, one can show that this equivalence does not depend on the choice of such $\xi$. 

However, we note here, that since $\Omega_0 G\subset P_0G$, the Hilbert space $H^{0}_p$ is a representation space for $\Omega_0 G$ too. Therefore we have a Proposition, arising from the direct integral decomposition of $H^{0}_p$.\\

\noindent {\bf Proposition 3.2.} {\it The representation $\rho_{\lambda}$ of $\Omega_0G$ on $H^0_p$ is the direct integral of the representations of $\Omega_0G$: 
\begin{eqnarray}
\rho_{\lambda}=\int^{\oplus}_{\mathbb{R}} \rho^{X}_{\lambda}dX, 
\end{eqnarray}
where $\rho^{X}_{\lambda}$ is equivalent to $\rho^0_{\lambda_X}$ such that  
$\lambda_X=e^{\xi}\lambda$, $\xi\in C'[0,2\pi]$ and $\xi(0)=0$, $\xi(2\pi)=X$.\\
}

\noindent In order to define the representation of $PG$ and $\Omega G$, one needs to extend the Hilbert space, i.e. one has to consider the Hilbert spaces $H_p=L^2(C_{0}[0,2\pi]\oplus \mathbb{R}, dw^t\times dx_0;\mathbb{C})$ and 
$H_l=L^2(C_{0,0}[0,2\pi]\oplus \mathbb{R}, dw^t\times dx_0;\mathbb{C})$, where $dx_0$ is the Lebesgue measure on $\mathbb{R}$.
We also note that every element $(e^{\alpha}, b)$ of $PG$ (resp. $\Omega G$) can be represented as 
$(e^{\alpha}, b)=(e^{\alpha_0}e^{\t \alpha}, b)$, where $\alpha_0=\alpha(0)$ and $\t \alpha(0)=0$.  
Then the following theorem is true.\\

\noindent {\bf Theorem 3.2.} {\it i) The following action of the group element $g(e^{\alpha}, b)$ on the space of 
functionals of $x(u), x_0$:
\begin{eqnarray}\label{act2}
&&\t{\rho^p}_{\lambda}(e^{\alpha},b)(f)(x, x_0)=\\
&&e^{-\frac{1}{4t}\int_0^{2\pi}\alpha'(u)\alpha'(u)\ud u-\frac{1}{2t} \int_0^{2\pi}\alpha'(u)\ud x(u)}
e^{\int_0^{2\pi} \lambda(u)b(u)e^{x(u)+x_0}\ud u}f(x+\t\alpha, x_0+\alpha_0)\nonumber
\end{eqnarray}
defines a unitary representation of $PG$  on the Hilbert space $H_p$ iff $\lambda\in iL^2([0,2\pi]$. \\
ii) The formula \rf{act2} defines the unitary representation 
$\t \rho^l_{\lambda}$ of $\Omega G$ on the Hilbert space $H_l$ iff $\lambda \in iL^2([0,2\pi];\mathbb{R})$.}\\

\noindent Note, that the representations $\rho^p_{\lambda}$,  $\rho^l_{\lambda}$ are equivalent to 
$\rho^p_{e^t\lambda}$,  $\rho^l_{e^t\lambda}$ correspondingly, for any $t\in \mathbb{R}$. 
In particular, in the case when $\lambda=const$ 
there exist only three inequivalent classes of unitary representations of $PG$ ($\Omega G$), 
corresponding to  $\t \rho_0$, $\t\rho_i$, $\t\rho_{-i}$ respectively, like it was in the finite dimensional case. \\

\noindent{\bf 3.3. Central extension and Lie algebra generators.} In this subsection, we will construct the unitary representations of groups 
$\hat{G}$ on $H_l$. Let $(e^{\alpha},b, s)$ be an element of $\hat G$. Let us consider the following action of this element on the function $H_l$:
\begin{eqnarray}\label{act3}
&&\rho^l_{\lambda,k}((e^{\alpha},b, s))(f)(x,x_0)=
e^{is}e^{-\frac{1}{4t}\int_0^{2\pi}\alpha'(u)\alpha'(u)\ud u-\frac{1}{2t} \int_0^{2\pi}\alpha'(u)\ud x(u)}\nonumber\\
&&e^{ik\int_0^{2\pi}\alpha(u)dx(u)}e^{\int_0^{2\pi} \lambda(u)b(u)e^{x(u)+x_0}\ud u}f(x+\t \alpha,x_0+\alpha_0).
\end{eqnarray}
{\bf Theorem 3.3.} {\it Formula \rf{act3} defines a unitary representation of $\hat G$ for the central charge $k\in \mathbb{R}$, and $\lambda\in i L^2([0,2\pi];\mathbb{R})$. The representation ${\rho^l}_{e^{\xi}\lambda,k}$ is equivalent to ${\rho^l}_{\lambda,k}$ if $\xi$ is any absolutely continuous function on $[0,2\pi]$, such that $\xi(0)=\xi(2\pi).$}\\ 

Ignoring $x_0$ dependence in \rf{act3}, one can define the unitary representation of $\hat{G}_0$ on $H^0_l$, for 
which we will keep the same notation $\rho^l_{\lambda,k}$. We also mention that the operators $\rho^l_{\lambda,k}(e^{\alpha_1}, 0, s_1)$ and $\rho^l_{\lambda,-k}(e^{\alpha_2}, 0, s_2)$  corresponding to the action of the $A$ subgroup of $\hat{G}$ commute.

One can write down the expressions for the Lie algebra generators. As we can see, the generator which is responsible for the action of the "$B$-subgroup"  has the form:
\begin{eqnarray}
T_bf(x,x_0)=\frac{d}{d\epsilon}_{\vert_{\epsilon=0}}\rho^l_{\lambda,k}(1,\epsilon b,0)=\int_0^{2\pi}\lambda(u)b(u)e^{x(u)+x_0}du.
\end{eqnarray}
It is well-defined for all the functionals $f(x,x_0)\in H_l$ such that $e^{x_0}f(x,x_0)\in H_l$. However, the $a$-subgroup generators are more peculiar:
\begin{eqnarray}
&&D_{k,\alpha}f(x)=\frac{d}{ds}_{\vert_{\epsilon=0}}\rho^l_{\lambda,k}(e^{\epsilon\alpha},0,0)=\\
&&(ik\int_0^{2\pi}\alpha(u) dx(u)+\frac{1}{2t}
\int_0^{2\pi}\alpha'(u)dx(u))f(x)+\frac{d}{d\epsilon}_{|_{\epsilon=0}}f(x+\epsilon\t \alpha, x_0+\epsilon \alpha_0).\nonumber
\end{eqnarray}
Therefore the action of the Lie algebra element $T_{\alpha}$ is defined only on the set of weakly differentiable functionals, which form a dense subset in $H_l$. In the notations of variational calculus, one can 
introduce the operators:
\begin{eqnarray}\label{var}
&&D_{\alpha,k}(\alpha)=\alpha_0\frac{\p}{\p x_0}+\int_0^{2\pi}\t\alpha(u)\frac{\delta}{\delta x(u)}du+
\frac{1}{2t}\int_0^{2\pi}\t\alpha'(u)x'(u)du+\nonumber\\
&& ik\int_0^{2\pi} du\t\alpha(u)x'(u)du,  \nonumber\\
&&T_b=\int_0^{2\pi}\lambda(u)b(u)e^{x_0+x(u)}du.
\end{eqnarray}
 One can see that \rf{var} satisfies the needed relation:
\begin{eqnarray}
[D_{\alpha_1,k}, D_{\alpha_2,k}]=-2ik \int_0^{2\pi}\alpha_1'(u)\alpha_2(u)du, 
\quad [D_{\alpha,k}, T_b]=T_{\alpha b}.
\end{eqnarray}
  
Finally, let us consider the following two semigroups of $\hat{G}_+,\hat{G}_-\in\hat G$, such that they 
consist of group elements $g(e^\alpha, b)$, where $b(u)>0$ or $b(u)<0$ for all $u$ correspondingly. The following statement will be important will be important for the construction of loop $\Gamma$-function in the next section.\\
 
\noindent{\bf Proposition 3.3.} {\it Let $Im\lambda, Re\lambda\in L^2([0,2\pi];\mathbb{R})$. Then \rf{act3} defines a representation of 
the semigroup $\hat{G}_+$ (resp. $\hat{G}_-$) if $Re \lambda(u)<0$ (resp.  $\lambda>0$) on $[0,2\pi]$.}\\
  
\noindent In the following we will denote the representations of $\hat{G}_+,\hat{G}_-$ the same way, namely 
$\rho^l_{\lambda,k}$.\\ 

\noindent{\bf 3.4. Irreducibility and classification of unitary representations.} In Section 2, we learned that the 
irreducible unitary representations of $G$ turn out to be equivalent to either $R_{\pm i}$ or $T_{\rho}$. In the loop case the following weaker result holds \cite{zeit}. \\

\noindent{\bf Theorem 3.4.} {\it Representations $\rho^l_{\lambda,k}$, such that  $\lambda(u)\equiv\lambda\in \mathbb{R}\backslash 0$ of $\hat{G}$ ($\hat{G}_0$) on $H_l$ ($H^0_l$), are irreducible.}\\

One can prove that if $i\lambda$ is either a positive or negative function from $L^2([0,2\pi],\mathbb{R})$, 
the representation $\rho^l_{\lambda,k}$ is also irreducible. 
However, if $\lambda(u)=0$  $\rho^l_{\lambda,k}$ reduces to the representation of the $A$-subgroup, i.e. the representations of the loop Heisenberg group. It is not irreducible since $[\rho^l_{\lambda,k}(e^{\alpha_1}, 0, s), \rho^l_{\lambda,-k}(e^{\alpha_2}, 0, t)]=0$. To author's knowledge the classification of the unitary representations of the loop Heisenberg group is not yet known (see e.g. \cite{ccr} for review of the subject). 
The same argument as in Theorem 3.4, shows that if $\lambda(u)=0$ for $u\in [a,b]$ the representation $\rho^l_{\lambda,k}$ is not irreducible. Therefore, we make the following conjecture about the classification of the irreducible representations 
of $\hat G$. 
Hence the irreducible unitary representations of $\hat G$ which we considered, are either equivalent to the  representations of the $A$-subgroup, or equivalent to the representations $\rho^l_{\lambda,k}$, where $i\lambda\in L^2([0,2\pi], \mathbb{R})$ is strictly positive or negative function on $[0,2\pi]$. Moreover, we know that  
$\rho^l_{\lambda,k}$ and $\rho^l_{e^{\xi}\lambda,k}$ are equivalent to each other if $\xi\in C'_{0,0}[0,2\pi]$. 
Therefore, an interesting problem to study is the classification of all finite-dimensional representations of $\hat{G}$. A reasonable conjecture (by analogy with the finite-dimensional case) would be that the three discussed classes of representations, namely $\rho^l_{\lambda,k}$ for positive and negative $\lambda$ and irreducible unitary representations of the loop Heisenberg group exhaust all irreducible unitary representations of $\hat{G}$. \\

\section{(Loop) $\Gamma$-function and the action of the affine loop group}

\noindent{\bf 4.1. Fourier transform for the classical Wiener measure.} In this section, we consider the generalizations of the formula relating the action of the group $G$ and the $\Gamma$-function (see Section 1). In particular, we introduce a new object, which we will refer to as the $loop$ $\Gamma$-$function$. In order to do that one needs to construct the generalization of the Fourier/Laplace transform in the case of the Wiener measure. We have already seen a  unitary transformation on the $L^2$ space for the abstract Wiener measure, called the Fourier-Wiener transform, 
but we will choose another transformation here.

Let us consider the following transformation on the Hilbert space $H^0_l$:
\begin{eqnarray}
\mathcal{F}f(p)=\int_{C_{0,0}[0,2\pi]} e^{i\int_0^{2\pi}p(u)x(u)du} f(x)dw_0^t(x),
\end{eqnarray}
where $p(u)\in C[0,2\pi]$ and $p(0)=p(2\pi)=0$. Unlike the usual Fourier transform for the Lebesgue measure on the real 
line, the transformation $\mathcal{F}$ is not a unitary operator.\\

\noindent{\bf Proposition 4.1.} {\it The operator  $\mathcal{F}$ is a compact normal operator on  
$H^0_l$ with no zero eigenvalues.}\\

\noindent{\bf Proof.} The general condition \cite{conway} for the general integral operator\\  
$K: L^2(X,d\mu(x))\to L^2(Y,d\nu(y))$, such that 
\begin{eqnarray}
Kf(y)=\int_X K(y,x)d\mu(x)
\end{eqnarray}
to be compact, is that $\int \int  |K(x,y)|^2d\mu(x)d\nu(y)<\infty$. In order to prove that it is normal, i.e. $\mathcal{F}\mathcal{F}^*=\mathcal{F}^*\mathcal{F}$, one just needs to write explicitly the resulting expressions and then use 
the Fubini theorem. To show that the operator $\mathcal{F}$ has no nonzero eigenvalues, one needs to use the fact that the exponentials of the form $e^{\int_0^{2\pi}{\alpha'(u)dx(u)}}$, where $\alpha(u)$ is absolutely continuous, form a dense subset in $H^0_l$ \cite{hida}.   

\hfill $\blacksquare$

\noindent Taking into account that $\mathcal{F}$ is a normal operator and using the polar decomposition theorem, one can decompose it 
as $\mathcal{F}=U_{\mathcal{F}}K$, where $U_{\mathcal{F}}$ is a unitary operator on $H^0_l$ and $K$ is a compact self-adjoint operator, such that $K=\sqrt{\mathcal{F}\mathcal{F}^*}$. 

It is obvious that the operators $\mathcal{F}$ and $U_{\mathcal F}$ can be continued to the space $H_l$. 
At the same time, we define another unitary operator $F$, which is a standard Fourier transform with respect to measure $dx_0$. One can show that $\mathcal{F}$ and $F$ commute. \\

\noindent{\bf 4.2. Loop $\Gamma$-function}.
Let us consider the  following expression:
\begin{eqnarray}
\rho^{l,\mathcal{F}}_{\lambda,k}(e^{\alpha},b,s)=U_{\mathcal{F}}\rho^l_{\lambda,k}(e^{\alpha},b,s)U^*_{\mathcal{F}}.
\end{eqnarray}
Since $U_{\mathcal{F}}$ is unitary, $\rho^{l,\mathcal{F}}_{\lambda,k}$ defines an equivalent representation 
of $\hat G$ on $H_l$. Because of the results of subsection 4.1, on the image  of $K$, one can rewrite it as follows:
\begin{eqnarray}
\rho^{l,\mathcal{F}}_{\lambda,k}(e^{\alpha},b,s)=K^{-1}\mathcal{F}\rho^l_{\lambda,k}(e^{\alpha},b,s)\mathcal{F}^*K^{-1}.
\end{eqnarray}
Here $K$ is a fixed self-adjoint operator, so we are interested in the object $\mathcal{F}\rho^l_{\lambda,k}(e^{\alpha},b)\mathcal{F}^*$. We consider the case when $(e^{\alpha},b)$ in $\hat{G}_{+}$ and $Im\lambda=0$, $Re \lambda(u)>0$. Let us write it down explicitly:
\begin{eqnarray}
&&\mathcal{F}\rho^l_{\lambda,k}(e^{\alpha},b,s)\mathcal{F}^*f(x,x_0)=\nonumber\\
&&e^{is}e^{-\frac{1}{4t}\int_0^{2\pi}\alpha'(u)\alpha'(u)du}\int_{C_{0,0}[0,2\pi]} e^{-i\int_0^{2\pi}(p(u)x(u))du}
e^{-\frac{1}{2t} \int_0^{2\pi}\alpha'(u)\ud p(u)}\nonumber\\
&&e^{ik\int_0^{2\pi}\alpha(u)dp(u)}e^{\int_0^{2\pi} \lambda(u)b(u)e^{p(u)+x_0}\ud u}\nonumber\\
&&\int_{C_{0,0}[0,2\pi]} e^{i\int_0^{2\pi}(p(u)+\t\alpha(u)y(u)du} f(y,x_0+\alpha_0)dw_0^t(y)dw_0^t(p).
\end{eqnarray}
Using the Fubini theorem, one can rewrite it as follows:
\begin{eqnarray}
&&\mathcal{F}\rho^l_{\lambda,k}(e^{\alpha},b,s)\mathcal{F}^*f(x,x_0)=\nonumber\\
&&\int_{C_{0,0}[0,2\pi]} \mathbb{K}^{\lambda,k}_{s,\alpha,b}(x-y, x_0) e^{i\int_0^{2\pi}\t\alpha(u)y(u)du}f(y, x+\alpha_0)dw_0^t(y),
\end{eqnarray}
where 
\begin{eqnarray}
&&\mathbb{K}^{\lambda,k}_{s,\alpha,b}(x-y,x_0)=e^{is}e^{-\frac{1}{4t}\int_0^{2\pi}\alpha'(u)\alpha'(u)}\\
&&
\int_{C_{0,0}[0,2\pi]} e^{i\int_0^{2\pi}p(u)(x(u)-y(u))du}e^{\int_0^{2\pi} \lambda(u)b(u)e^{p(u)+x_0}\ud u}\nonumber\\
&&e^{ik\int_0^{2\pi}\alpha(u)dp(u)}e^{-\frac{1}{2t} \int_0^{2\pi}\alpha'(u)\ud p(u)}
dw_0^t(p).\nonumber
\end{eqnarray}
In the case when $\alpha$ is twice differentiable, one can see that the object $\mathbb{K}^{\lambda,k}_{s,\alpha,b}(z,x_0)$ up to factors independent of $x,y$ is a particular case of the following functional:
\begin{eqnarray}
\hat\Gamma_{\mu}(z)=\int_{C_{0,0}[0,2\pi]} e^{\int_0^{2\pi}p(u)z(u)du}e^{-\int_0^{2\pi} \mu(u)e^{p(u)}\ud u}dw_0^t(p),
\end{eqnarray}
where $Re z, Im z \in L^2([0,2\pi];\mathbb{R})$, $\mu\in L^2([0,2\pi];\mathbb{R})$. We will call $\hat\Gamma_{\mu}(z)$ the {\it loop} {\it Gamma} $function$ or simply $\hat \Gamma$-$functional$. It has the following properties \cite{zeit}.\\

\noindent {\bf Theorem 4.1.} \\
\noindent{\it i) $\hat\Gamma_{\mu}(z)$ is well defined for any $Re z, Im z \in L^2([0,2\pi];\mathbb{R})$, $\mu\in L^2([0,2\pi];\mathbb{R})$ and $\mu(u)\ge 0$ on $[0,2\pi]$. \\
ii) The following relation is valid: 
\begin{eqnarray}\label{gamma}
&&\int_0^{2\pi} g(v)\mu(v)\hat\Gamma_{\mu}(z+\delta_v)dv=\int _{0}^{2\pi}g(v)z(v)dv\hat\Gamma_{\mu}(z)+\nonumber\\
&&\frac{1}{t}\int_0^{2\pi}g''(v)\frac{\delta}{\delta z(v)}\hat\Gamma_{\mu}(z)dv, 
\end{eqnarray}
where $g(v)$ is any twice differentiable function on $[0,2\pi]$, such that $g(0)=g(2\pi)=0$, 
$\delta_v=\delta(u-v)$ is a delta-function on the interval $[0,2\pi]$ and \\
$
\int_0^{2\pi} \xi(v)\frac{\delta}{\delta z(v)}\hat\Gamma_{\mu}(z)=\frac{d}{d\epsilon}_{|\epsilon=0}\hat\Gamma_{\mu}(z+\epsilon\xi)$ for any continuous function $\xi$. }\\

\noindent{\bf Proof.} An important step in the proof is the 
consideration of the infinitesimal form of the translation invariance. Let $f\in H^0_l$ such that it is weakly differentiable. Then we have the following property:
\begin{eqnarray}
&&\int_{C_{0,0}[0,2\pi]} e^{-\frac{1}{2t}\int_0^{2\pi}\epsilon^2g'(u)g'(u)\ud u+\frac{1}{t} \int_0^{2\pi}\epsilon g''(u)x(u)\ud u}f(x+\epsilon g )dw_0^t(x)=\nonumber\\
&&\int f(x)  dw^t(x),
\end{eqnarray}
where $g\in C_{0,0}^2[0,2\pi]$.
This is just the translation property (the translation is with respect to function $\epsilon g$), where $\epsilon$ is some real parameter.  Then, if we differentiate with respect to $\epsilon$ at zero, we obtain the following formula:
\begin{eqnarray}\label{rel}
&&\int_{C_{0,0}[0,2\pi]}\Big(\frac{1}{t} \int_0^{2\pi}g''(u)x(u)\ud u\Big)f(x)dw_0^t(x)=\nonumber\\
&&-\int_{C_{0,0}[0,2\pi]} \frac{d}{d\epsilon}_{|_{\epsilon=0}}f(x+\epsilon g)dw_0^t(x).
\end{eqnarray}
Let us apply this property to the integrand of the $\hat \Gamma$-functional, i.e. let 
$f(x)=F(x,z)\equiv e^{\int_0^{2\pi}x(u)z(u)du}e^{-\int_0^{2\pi} \mu(u)e^{x(u)}\ud u}$. Then we have:
\begin{eqnarray}\label{rel1}
&&\int_{C_{0,0}[0,2\pi]}\Big(\frac{1}{t} \int_0^{2\pi}g''(u)x(u)\ud u\Big)F(x,z)dw_0^t(x)=\nonumber\\
&&\frac{1}{t}\int_0^{2\pi}g''(v)\frac{\delta}{\delta z(v)}\hat\Gamma_{\mu}(z),
\end{eqnarray}
and
\begin{eqnarray}\label{rel2}
&&-\int_{C_{0,0}[0,2\pi]} \frac{d}{d\epsilon}_{|_{\epsilon=0}}F(x+\epsilon g,z)dw_0^t(x)=\nonumber\\
&&-\int_{C_{0,0}[0,2\pi]}\Big(\int^{2\pi}_0 g(v)z(v)dv-\int^{2\pi}_0\mu(v)g(v)e^{x(v)}dv\Big)F(x,z)dw_0^t(x)=\nonumber\\
&&-\int^{2\pi}_0 g(v)z(v)dv\hat\Gamma_{\mu}(z)+\int^{2\pi}_0g(v)\mu(v)\hat{\Gamma}(z+\delta_v)dv.
\end{eqnarray}
Therefore, combining \rf{rel}, \rf{rel1} and \rf{rel2} we obtain \rf{gamma}. 
\hfill $\blacksquare$\\

The property ii) from the theorem above is a natural generalization of the property of the ordinary $\Gamma$-function: 
$\Gamma(z+1)=z\Gamma(z)$. We notice, however, that there is also an extra term, depending on the $t$-parameter, related to the choice of the measure. So, the proper analogue of the functional $\hat{\Gamma}$ would be the "regularized" $\Gamma$-function:
\begin{eqnarray}
\Gamma_{\mu,t}(z)=\int_{\mathbb{R}} e^{-\mu e^x} e^{zx}e^{-\frac{x^2}{2t}}dx,
\end{eqnarray}
where $\mu, t\in{R_{+}}$. It satisfies the equation
\begin{eqnarray}
\mu\Gamma_{\mu,t}(z+1)=z\Gamma_{\mu,t}(z)-\frac{1}{t}\frac{d\Gamma(z)}{dz}.
\end{eqnarray}
This function is well defined for all complex values of  $z$. Its relation to the original $\Gamma$-function can be obtained by considering $t\to \infty$ limit: 
\begin{eqnarray}
\lim_{t\to \infty}\Gamma_{\mu,t}(z+1)=\mu^{-z}\Gamma(z).
\end{eqnarray}
Also, we note here that this function can be related to the matrix elements of the representations of $G$ (see Section 2) if we would consider the Gaussian measure instead of the Lebesgue measure on the real line. 
\\

\section{Representations of $\widehat{\mathcal{A}}$ and $\widehat{\mathcal{K}}$}
\noindent {\bf 5.1. Definitions.} In this section we consider the following *-algebras, which we denote $\widehat{\mathcal{A}}$ and $\widehat{\mathcal{K}}$ that are the affine analogues of $\mathcal{A}, \mathcal{K}$, considered in Section 2. The first algebra is close to the Lie algebra of a loop $ax+b$-group (see Section 3), which is generated by $h_n$, $e^{\pm}_n$, $n\in \mathbb{Z}$,
so that the generating "currents" are:
\begin{eqnarray}
h(u)=\sum_{n\in\mathbb{Z}} h_{-n} e^{inu}, \quad e^{\pm}(u)=\sum_{n\in\mathbb{Z}} e^{\pm}_{-n}e^{inu} 
\end{eqnarray}
and obey the following commutation relations, expressed via the generating functions:
\begin{eqnarray}\label{relcom}
&&[h(u), h(v)]=0, \quad [e^{\pm}(u), e^{\pm}(v)]=0, \quad [e^{\pm}(u), e^{\mp}(v)]=0,\nonumber\\
&&[h(u), e^{\pm}(v)]=\pm i\delta(u-v)e^{\pm}(v), \quad e^{\pm}(u)\cdot e^{\mp}(u)=1.
\end{eqnarray}
The *-structure is such that $h(u)$, $e^{\pm}(u)$ are Hermitian, i.e.
\begin{eqnarray}
h(u)^*=h(u),\quad e^{\pm}(u)^*=e^{\pm}(u).
\end{eqnarray}
The algebra $\widehat{\mathcal{K}}$ is generated by $h_n$, $\alpha^{\pm}_n$, $n\in \mathbb{Z}$,
it has similar commutation relations, but a different *-structure:
\begin{eqnarray}
&&h(u)=\sum_n h_{-n} e^{inu}, \quad \alpha^{\pm}(u)=\sum_n \alpha^{\pm}_{-n}e^{inu},\nonumber\\ 
&&[h(u), h(v)]=0, \quad [\alpha^{\pm}(u), \alpha^{\pm}(v)]=0, \quad [\alpha^{\pm}(u), \alpha^{\mp}(v)]=0,\nonumber\\
&&[h(u), \alpha^{\pm}(v)]=\mp \delta(u-v)\alpha^{\pm}(v), \quad \alpha^{\pm}(u) \alpha^{\mp}(u)=1,\nonumber\\
&&h(u)^*=h(u),\quad \alpha^{\pm}(u)^*=\alpha^{\mp}(u).
\end{eqnarray}
Using Gaussian integration on Hilbert spaces, in the next subsection we will construct some representations of these algebras.\\

\noindent{\bf 5.2. Construction of representations.} In order to construct the representations of $\widehat{\mathcal{A}}$ and $\widehat{\mathcal{K}}$ one can use the Gaussian measure on a Hilbert space (see Appendix). 
Let us consider the Fourier series of a function from $L^2(S^1, \mathbb{R})$:
\begin{eqnarray}
x(u)=\sum_{n\in \mathbb{Z}}x_{-n}e^{inu}, \quad x_0\in \mathbb{R}, \quad x_n^*=x_{-n}.
\end{eqnarray}
Let us introduce two quadratic forms defining two types of trace-class operators on $L^2(S^1, \mathbb{R})$, which will determine the appropriate Gaussian measures:
\begin{eqnarray}\label{bil}
&& B_A(x,x)=\frac{1}{2}\sum_{n\ge 1}\xi_n^{-1}x_nx_{-n}+\xi_0^{-1}x_0^2,\nonumber\\
&& B_{K}(x,x)=\frac{1}{2}\sum_{n\ge 1}\xi_n^{-1}x_nx_{-n},
\end{eqnarray}
where $\xi_n>0$ for all $n$ and $\sum_n\xi_n<\infty$. 
The Gaussian measures we are interested in, heuristically can be expressed as follows:
\begin{eqnarray}
&& dw_{A}=(\sqrt{det(2\pi N_A)})^{-1}e^{-B_A(x,x)}dx_0\prod^{\infty}_{n=1}[\frac{i}{2}dx_n\wedge dx_{-n}],\nonumber \\
&& dw_{K}=(\sqrt{det(2 \pi N_{K})})^{-1}e^{-B_K(x,x)}d\phi\prod^{\infty}_{n=1}[\frac{i}{2}dx_n\wedge dx_{-n}],
\end{eqnarray}
where $N_A, N_{K}$  are trace-class diagonal operators determined by the quadratics forms \rf{bil}. Here as before the range for $\phi$ is $[0,2\pi]$. Literally the difference between two measures is that in the second one we compactified the zero mode $x_0$ on a circle with parameter $\phi$ for $dw_{K}$. Hilbert spaces of square-integrable functions with respect to these measures are denoted in the following as $\mathcal{H}_A$ and $\mathcal{H}_K$. 

 Let us construct the unitary representations of $\widehat{\mathcal{A}}$ and $\widehat{\mathcal{K}}$-algebras in the Hilbert spaces $\mathcal{H}_A$ and $\mathcal{H}_K$, correspondingly. 
 
 We will explicitly define the operators, which will represent the generators. Let us start from  $\widehat{\mathcal{A}}$-algebra. First of all, let us extend the index of $\xi_n$ to all integers, so that $\xi_{-n}=\xi_n$, $n\in \mathbb{Z}$. We need the following differential operators: 
 \begin{eqnarray}\label{ab}
b_{-n}=i(\p_n-\xi_n^{-1}x_{-n}), 
\quad a_{-n}=i\p_n,
\end{eqnarray}
where $\p_{n}=\frac{\p}{\p x_n}$. These operators are formally conjugate to each other, 
\begin{eqnarray}
a_n^*=b_{-n},
\end{eqnarray}
being considered on a certain dense set $D^l_A$, i.e. functions which are the sums of monomials 
\begin{eqnarray}
\prod^n_{k=1}\langle\mu_k,x\rangle\prod^m_{s=1} \langle \lambda_s, e^{\pm x}\rangle, 
\end{eqnarray}
where $\langle,\rangle$ is the standard $L^2(S^1, \mathbb{R})$ pairing and $\mu_k, \lambda_k$ are trigonometric polynomials. 
We define the operators $h_n$ as follows:
\begin{eqnarray}\label{ha}
h_{-n}=\frac{1}{2}(a_{-n}+b_{-n})=i(\p_n-\frac{1}{2}\xi_n^{-1}x_{-n}).
\end{eqnarray}
Hence, on a dense set $h_n^*=h_{-n}$, so that the current $h(u)=\sum_{n\in \mathbb{Z}}h_{-n}e^{inu}$
is Hermitian. 
Also, note that Hermitian currents $h(u), x(v)$ generate the infinite-dimensional Heisenberg algebra:
\begin{eqnarray}
[h(u), x(v)]=i\delta(u-v).
\end{eqnarray}
We also define the currents $e^{\pm}(u)$: 
\begin{eqnarray}\label{ea}
e^{\pm}(u)=e^{\pm x(u)}.
\end{eqnarray}
It also follows that they satisfy the commutation relations \rf{relcom}. One can show that $e^{\pm x(u)}$ for any $u$ are Hermitian operators, considered on the dense set $D^l_A$. Therefore, this gives a unitary representation of $\widehat{\mathcal{A}}$-algebra. 

Similarly one can construct the unitary representations of $\widehat{\mathcal{K}}$-algebra. Let us consider the dense set $D^l_A$ in $\mathcal{H}_K$ of the following form:
\begin{eqnarray}
\prod^n_{k=1}\langle\mu_k,x\rangle\prod^m_{s=1} \langle \lambda_s, e^{\pm ix^c}\rangle, 
\end{eqnarray}
where $\lambda_i$ are trigonometric polynomials without the constant term and $x^c(u)=\phi+\sum_{n\neq 0}x_n e^{-inu}$.
The currents $h(u), \alpha^{\pm}(u)$ are defined by the following formulas:
\begin{eqnarray}\label{ag}
&& h(u)=\sum_{n\neq 0} i(\p_n-\frac{1}{2}\xi_n^{-1}x_{-n})e^{inu}+i\p_{\phi},\nonumber\\
&& \alpha^{\pm}(u)=e^{\pm ix^c(u)}.
\end{eqnarray}
It is possible to introduce operators $a_n, b_n$ and  define them by
the same formulas as in \rf{ab} for all $n\neq 0$. For $n=0$ we put $a_0=b_0=i\p_{\phi}$. 
Therefore, the elements of the form $a\cdot v_0$, where $v_0\equiv 1\in \mathcal{H}_K$ and $a$ belongs to universal enveloping algebra of  $\widehat{\mathcal{K}}$ and the action of 
the generators is given by the formulas \rf{ag}. Such elements generate a dense set 
$D^l_K$ in $\mathcal{H}_K$, which is a unitary representation of  $\widehat{\mathcal{K}}$-algebra.

An important notion which is necessary for our construction is the correlator associated with the representation. 

 By the correlator of generators $T_1, ..., T_n$ of $\widehat{\mathcal{A}}$-algebra (resp. $\widehat{\mathcal{K}}$-algebra) 
 we mean the following expression:
 \begin{eqnarray}
 <T_1...T_n>\equiv \langle v_0, T_1...T_n v_0 \rangle,
 \end{eqnarray}
 where the pairing $\langle, \rangle $ is of the Hilbert space $\mathcal{H}_A$ (resp. $\mathcal{H}_K$), $v_0$ is the vector corresponding to the constant function $1$ and  $T_1, ..., T_n$ 
 are the generators $h_n, e^{\pm}_m$ (resp. $h_n, \alpha^{\pm}_m$). 
 
We remind that $h_n=\frac{1}{2}(a_n+b_n)$.  It is crucial for the  following, that the correlators  
\begin{eqnarray}\label{van}
<T_1...T_na_k>\quad , \quad 
<b_kT_1...T_n>
\end{eqnarray}
vanish for any generators $T_1, ...,T_n$, as a simple consequence of properties of the Gaussian integration.

It is natural to call operators $a_n$  $annihilation$ operators and 
$b_n$  $creation$ operators. This allows us to define the normal ordering. Namely, when we write down the expression $:T_1...T_n:$ for the product of $n$ generators, we reorder them in such a way that 
the creation operators will be to the left and the annihilation ones to the right. 

This procedure together with the vanishing of the correlators \rf{van} also gives an easy method to compute the correlators: by means of commutation relations of generators, one can reduce the products $T_1, ...,T_n$ to the normally ordered expressions. Therefore, the result  will reduce to the correlators of the generators $e^{\pm}_n$ or $\alpha^{\pm}_n$. 
 We also note that in the case of $\widehat{\mathcal{K}}$-algebra the corresponding correlators are nonzero only if they have an equal number of $\alpha^+$ and $\alpha^-$ generators. 
\\
In order to compute the correlators of these generators it is easier to consider the appropriate currents instead of modes and use the Gaussian integration.

There are the following expressions for correlation functions:    
\begin{eqnarray}\label{core}
&&\langle e_+(u_1)...e_+(u_n)e_-(v_1)...e_-(v_m)\rangle=\nonumber\\
&&\exp(\sum^n_{i<j; i, j=1}N_A(u_i,u_j)+\sum^m_{r<s; r, s=1}N_A(v_r,v_s)-\sum^n_{k=1}\sum^m_{l=1}N_A(u_k,v_l)+\nonumber\\
&&\frac{n+m}{2}N_A(0,0)),\\
&&\label{corealpha}\langle \alpha_+(u_1)...\alpha_+(u_n)\alpha_-(v_1)...\alpha_-(v_m)\rangle=\nonumber\\
&&\delta_{n, m}\exp(-\sum^n_{i<j; i, j=1}N_K(u_i,u_j)-\sum^n_{i<j; i, j=1}N_K(v_i,v_j)+\sum^n_{k,l=1}N_K(u_k,v_l)\nonumber\\
&&+nN_K(0,0)),
\end{eqnarray}
where 
\begin{eqnarray}
&&N_A(u, v)=2\sum_{n\ge 0}cos(n(u-v))\xi_n, \nonumber\\ 
&&N_K(u, v)=2\sum_{n>0}cos(n(u-v))\xi_n.
\end{eqnarray}

We note here that in the case of correlators of the currents involving not only $e^{\pm}, \alpha^{\pm}$, but 
also $h(u)$, the resulting expression will consist of monomials \rf{corealpha} multiplied on a certain product of delta-functions, coming from commutation relations of $\widehat{\mathcal{K}}$, $\widehat{\mathcal{A}}$-algebras.

We have shown above, that using the Gaussian measure, one can construct unitary representations of $\widehat{\mathcal{A}}$, 
$\widehat{\mathcal{K}}$-algebras.  However, it is possible to simplify those representations by making $a_k$ commute with $b_s$ for any $k$ and $s$. However, the resulting module will not be unitary, i.e. the pairing though nondegenerate will lose its positivity. 

At first, let us give the explicit description of such module for $\widehat{\mathcal{K}}$-algebra. For this we consider the vacuum 
vector $v_0$, such that $v_0$ is annihilated by $a_k$:
\begin{eqnarray}
a_kv_0=0.
\end{eqnarray}
Then the module which we will refer to as $\mathcal{V}_A$ is spanned by the following vectors:
\begin{eqnarray}
b_{m_1}...b_{m_s}e^{\pm}_{n_1}...e^{\pm}_{n_r}v_0,
\end{eqnarray}
 where $n_1, ..., n_r, m_1, ..., m_s\in \mathbb{Z}$. 
 Let us  begin to define the pairing with the postulation of
the points:
$e(u)=\sum_ne_ne^{inu}$ is a Hermitian current, $b(u)^*=a(u)$, so that $b(u)=\sum_n b_ne^{-inu}$, $a(u)=\sum_n a_ne^{-inu}$ and 
\begin{eqnarray}
[b(u), a(v)]=0.
\end{eqnarray}
The pairing is uniquely defined by the correlator of $e^{\pm}$ currents is given by \rf{core}. 
Similarly one can define the module $\mathcal{V}_K$ for $\widehat{\mathcal{K}}$ algebra
with the same conditions, just replacing $e^{\pm}$ with $\alpha^{\pm}$, as a result $\mathcal{V}_K$ is spanned by
\begin{eqnarray}
b_{m_1}...b_{m_s}\alpha^{\pm}_{n_1}...\alpha^{\pm}_{n_r}v_0.
\end{eqnarray}
One can define the pairing on $\mathcal{V}_K$ which is uniquely determined by the correlators of $\alpha^{\pm}$ currents \rf{corealpha}. 
  This pairing is Hermitean and nondegenerate. It gives a structure of nonunitary representations of the *-algebras $\widehat{\mathcal{A}}$, $\widehat{\mathcal{K}}$ on the spaces 
$\mathcal{V}_A$, $\mathcal{V}_K$  correspondingly. 

We note here that clearly, these representations are nonunitary, because you can easily find vectors $v\in \mathcal{V}_A, \mathcal{V}_K$, such that $\langle v,v\rangle=0$.

\section{Construction of representations for $\widehat{sl(2,\mathbb{R})}$}
\noindent {\bf 6.1. Regularization and commutator.} In order to construct the $\widehat{sl(2,\mathbb{R})}$ representations via $\widehat{\mathcal{A}}$, $\widehat{\mathcal{K}}$, one has to find an affine analogue of the formulas for $E,F, H$ and $J^3, J^{\pm}$. In the case of nontrivial central extension, the $\widehat{\mathcal{A}}$, $\widehat{\mathcal{K}}$ representations appear to be insufficient, one has to introduce the representation of the infinite-dimensional Heisenberg algebra, so that the generating current is $\rho(u)=\sum_{n\in \mathbb{Z}}\rho_ne^{-inu}$ and the commutation relations are:
\begin{eqnarray}
[\rho_n, \rho_m]=2\kappa n\delta_{n,-m},
\end{eqnarray}
where $\kappa\in \mathbb{R}_{>0}$. 
The irreducible module, the so-called Fock module $F_{\kappa,p}$ of this algebra is defined as follows. We introduce a vector $vac_p$ with the property $\rho_n vac_p=0$, $p\in \mathbb{R}$, $n>0$ so that
\begin{eqnarray}
F_{\kappa,p}=\{\rho_{-n_1}...\rho_{-n_k}vac_p; \quad n_1,..., n_k>0,\quad  \rho_0 vac_p=p\cdot vac_p\}.
\end{eqnarray}
The Hermitian pairing is defined so that 
\begin{eqnarray}
\langle vac_p, vac_p \rangle=1, \quad \rho^*_n=\rho_{-n}.
\end{eqnarray}
Another object required in this section is the regularized version of the 
$\rho$, $h$, 
$e^{\pm}$, 
$\alpha^{\pm}$ currents. Namely, for any $\varphi$, which stands for any of $\rho$, $h$, 
$x(u)$ or $x^c(u)$ we consider 
\begin{eqnarray}
\varphi(z,\bar z)=\sum_{n\ge 0}\varphi_n \bar z^n+\sum_{n> 0}\varphi_{-n} z^n,
\end{eqnarray}
where $z=re^{iu}$, so that $0<r\le 1$. We denote $e^{\pm}(z,\b z)\equiv e^{\pm x(z,\b z)}$ and $\alpha^{\pm}(z,\b z)\equiv e^{\pm i x^c(z,\b z)}$. 
The Wick theorem implies that the correlators of the regularized Heisenberg currents $\rho(z_1, \bar z_1), ...., \rho(z_n,\bar z_n)$ are finite as long as $0<|z_i|< 1$, i.e. 
the expressions
\begin{eqnarray}\label{heicor}
\langle vac_p, \rho(z_1,\bar z_1)...\rho(z_n,\bar z_n) vac_p\rangle
\end{eqnarray}
are finite. One can consider the limit  $|z_i|\to1$ in the sense of distributions, so that \rf{heicor} is the sum of products of distributions, i.e. 
\begin{eqnarray}\label{2pointf}
\langle vac_p, \rho(u_1)\rho(u_2)vac_p\rangle=\frac{2\kappa}{(1-e^{i(u_2-u_1-i0)})^2}+p^2.
\end{eqnarray}
Next we consider the following sets of composite regularized currents 
\begin{eqnarray}\label{sl2R1}
&&J^{\pm}(z,\bar z)=\\
&&\frac{i}{2}(b(z,\bar z)\alpha^{\pm}(z,\bar z)+\alpha^{\pm}(z,\bar z)a(z,\bar z))\pm\kappa\p_u\alpha^{\pm}(z,\bar z)
\pm \rho(z,\bar z)\alpha^{\pm}(z,\bar z),\nonumber\\
&&J^3(z,\bar z)=2i h(z,\bar z)-2\kappa\alpha^-(z,\bar z)\p_u \alpha^+(z,\bar z),\nonumber
\end{eqnarray}
and 
\begin{eqnarray}\label{sl2R2}
&&E(z,\bar z)=\\
&&\frac{i}{2}(b(z,\bar z)e^{+}(z,\bar z)+e^{+}(z,\bar z)a(z,\bar z))+i\kappa\p_ue^{+}(z,\bar z)
+ i\rho(z,\bar z)e^{+}(z,\bar z),\nonumber\\
&&F(z,\bar z)=\nonumber\\
&&-\frac{i}{2}(b(z,\bar z)e^{-}(z,\bar z)+e^{-}(z,\bar z)a(z,\bar z))+i\kappa\p_ue^{-}(z,\bar z)
+i \rho(z,\bar z)e^{-}(z,\bar z),\nonumber\\
&&H(z,\bar z)=-2i h(z,\bar z)+2i\kappa e^-(z,\bar z)\p_u e^+(z,\bar z).\nonumber
\end{eqnarray}
satisfying the Hermicity conditions:
\begin{eqnarray}
&&E(z,\bar z)^*=-E(z,\bar z), \quad F(z,\bar z)^*=-F(z,\bar z), \quad H(z,\bar z)^*=-H(z,\bar z),\nonumber\\
&&J^3(z,\bar z)^*=-J^3(z,\bar z), \quad J^{\pm}(z,\bar z)^*=-J^{\mp}(z,\bar z).
\end{eqnarray}
Then if  $\phi_k$ denote $E, F, H$ or $J^3, J^{\pm}$ then the correlators
\begin{eqnarray}\label{gencor}
\quad \langle \phi_{1}(z_1,\bar z_1)...\phi_{n}(z_n,\bar z_n)\rangle_p\equiv\langle v_0\otimes vac_p,\phi_{1}(z_1,\bar z_1)...\phi_{n}(z_n,\bar z_n)v_0\otimes vac_p\rangle,
\end{eqnarray}
 are well-defined for $0<|z_i|<1$. Actually, after the normal ordering procedure $\phi_{1}(z_1,\bar z_1)...\phi_{n}(z_n,\bar z_n)$ will be represented as normally ordered products with coefficients which are continuous functions of $u_1,..., u_n$, where 
$z_i=r_ie^{u_i}$. Therefore, the expression $\langle \phi_{1}(z_1,\bar z_1)...\phi_{n}(z_n,\bar z_n)\rangle_p$ will be given by the sum of correlators of exponentials $\alpha^{\pm}$ or $e^{\pm}$ which are well-defined even on the unit circle as we know from the formulas \rf{core}. 
Notice, that if only one of the 
currents $\phi_k(z_k, \bar z_k)$ is such that 
$|z_k|=e^{iu_k}$, i.e. it lies on the unit circle, the situation will not change. This is because during the normal ordering procedure, the commutators of $a$- and $b$- parts of this 
generator with other ones again produce continuous functions of $u_1,..., u_n$, since all other $z_i$ are 
inside the unit circle and the normal ordering wouldn't produce any distributions.
 
This allows us to make sense of the commutator of two currents on the unit circle. 
We consider the following limit of the difference of two correlators 
\begin{eqnarray}\label{comm}
&&\langle \phi_{1}(z_1,\bar z_1)...[\xi(u_1),\eta(u_2)]...\phi_{n}(z_n,\bar z_n)\rangle_p\equiv\\
&&\lim_{r_1,r_2\to 1}\big(\langle \phi_{1}(z_1,\bar z_1)...\xi(w_1,\bar w_1)\eta(w_2,\bar w_2)...\phi_{n}(z_n,\bar z_n)\rangle_p-\nonumber\\
&&\langle \phi_{1}(z_1,\bar z_1)...\eta(w_2,\bar w_2)\xi(w_1,\bar w_1)...\phi_{n}(z_n,\bar z_n)\rangle_p\big ),\nonumber 
\end{eqnarray} 
Then the following result holds (see \cite{fzsl2} Proposition 4.2, Theorem 4.1).
\\

\noindent {\bf Theorem 6.1.} {\it 
i) The expression \rf{comm} exists in the sense of distributions, more specifically the answer will contain a linear combination of delta-functions $\delta(u_1-u_2)$ and their derivatives. 
Here  $\phi_k$ stands for $E, F, H$ or $J^3, J^{\pm}$, $w_i=r_ie^{iu_i}$, $\xi(u_1)\equiv \xi(e^{iu_1},e^{-iu_1})$, $\eta(u_2)\equiv \eta(e^{iu_2},e^{-iu_2})$,\\ 
ii) The commutator \rf{comm} of the currents $\xi, \eta=E, F, H$ or $\xi,\eta=J^3, J^{\pm}$ exists and satisfies  the commutation relations for 
$\widehat{sl(2,\mathbb{R})}$ algebra with the central charge $\kappa$:
\begin{eqnarray}
&&[E(u),F(v)]=H(v)\delta(u-v)-4i\kappa\delta'(u-v), \quad [H(u), H(v)]=8i\kappa\delta'(u-v),\nonumber\\
&& [H(u), E(v)]=2E(v)\delta(u-v),\quad [H(u), F(v)]=-2F(v)\delta(u-v)
\end{eqnarray}
and
\begin{eqnarray}\label{sl2R}
&&[J^+(u),J^-(v)]=i J^3(v)\delta(u-v)+4i\kappa\delta'(u-v), \nonumber\\
&& [J^3(u), J^3(v)]=-8i\kappa\delta'(u-v),\nonumber\\
&&[J^3(u), J^{\pm}(v)]=\pm 2iJ^{\pm}(v)\delta(u-v).
\end{eqnarray}
}\\

It is important to mention that though we managed to define the commutator, based on the regularized commutators, this definition doesn't provide a representation, since the correlation functions of $E, F, H$ or $J^3, J^{\pm}$ do not exist, if more than one of the arguments lies on a circle. Then it is clear that the spaces from which we started, i.e. $\mathcal{H}_K\otimes F_p$ or  $\mathcal{H}_A\otimes F_p$ are not suitable to be spaces for $\widehat{sl(2,\mathbb{R})}$-module. However, we 
still have the regularized correlators which obey commutation relations. If we manage to eliminate divergencies in such a way that commutation relations and Hermicity conditions would be preserved, then the correlators will determine representation with the Hermitian bilinear form.  
In the next subsection we show a method, how to get rid of the  divergencies and redefine the correlator, so that it is well-defined when all the arguments are on a circle.\\

\noindent {\bf 6.2. Renormalization of correlators and construction of representations of $\widehat{sl}(2,\mathbb{R})$.} 
In this subsection, we will renormalize the correlators with currents on a circle for the cases of $\mathcal{V}_A$ and $\mathcal{V}_K$ representations. In both cases the description is very similar, so we focus on $\mathcal{V}_K$ case and the correlators of 
$J^3, J^{\pm}$ currents. Generalization to $\mathcal{V}_A$ and $E,F, H$ goes along the same path. 

In order to renormalize, at first we have to understand what kind of divergencies we are dealing with. For that purpose it is convenient to write down the expression for the correlator in the graphic form 
using the  
Feynman-like diagrams. 

Recall, that the easiest way to compute the correlator for $\widehat{\mathcal{K}}$ is to reduce the whole expression to the normal ordered form, i.e. all the creation operators $b(z,\bar z)$ are moved to the left and all the annihilation operators 
$a(z,\bar z)$ to the right.

Once we move the creation operator $b(z,\bar z)$, which is part of a  certain generator at the point to the left (or the annihilation operator $a(z,\bar z)$ to the right), it may produce the following terms arising from the commutation with the 
$\alpha^{\pm}$-generator of $J^{\pm}$ currents on the right (on the left in the case of $a(z,\bar z)$) of the given generator:
\begin{eqnarray}\label{lines}
&&[a(z,\bar z), \alpha^{\pm}(w,\bar w)]=\mp \alpha^{\pm}(w,\bar w)\delta(z,w),\nonumber\\
&&[\alpha^{\pm}(w,\bar w), b(z,\bar z)]= \pm \alpha^{\pm}(z,\bar z)\delta(z,w),
\end{eqnarray}
where $\delta(z,w)=\sum_{n\ge 0} (z\bar w)^n+\sum_{n> 0} (\bar zw)^n$. If $z,w$ are on the circle, i.e. $z=e^{iu}, w=e^{iv}$, we get $\delta(z,w)=\delta(u-v)$.
We will depict every term of the form \rf{lines}, which we obtain during the normal ordering procedure, as a line from one vertex to another: 

\begin{eqnarray}
\begin{xy}
(20,0)*+{\bullet}="2";
(20,5)*+{}="0";
(30,5)*+{\delta(z,w)}="1";
(40,0)*+{\bullet}="3";
"2";"3"  **\dir{-}; ?(.65)*\dir{>};
\end{xy}
\end{eqnarray}
Here the initial vertex corresponds to the term containing creation/annihilation operator and the terminal vertex correspond to the term containing $\alpha^{\pm}$. The direction of the arrow (to the right or to the left) indicates whether it was annihilation or creation operator.

At the same time, each generator $J^{\pm}$ contributes terms like that
\begin{eqnarray}
\alpha^{\pm}(z,\bar z)a(z,\bar z) , \quad b(z,\bar z)\alpha^{\pm}(z,\bar z).
\end{eqnarray}
They enter the diagram as vertices with one outgoing line to the right or to the left and any amount of incoming lines: 
\begin{eqnarray}
\begin{xy}
(0,0)*+{}="1";
(20,0)*+{\bullet}="2";
(20,3)*+{-}="-";
(20,5)*+{+}="+";
(40,0)*+{}="3";
(3,6)*+{}="1t";
(10,12)*+{}="1tt";
(3,-6)*+{}="1b";
(10,-12)*+{}="1bb";
 "1";"2" **\dir{-}; ?(.65)*\dir{>};
{\ar@{->} "2";"3"};
"1t";"2"  **\dir{-}; ?(.65)*\dir{>};
"1tt";"2"  **\dir{-}; ?(.65)*\dir{>};
"1b";"2"  **\dir{-}; ?(.65)*\dir{>};
"1bb";"2"  **\dir{-}; ?(.65)*\dir{>};
(60,0)*+{}="4";
(80,0)*+{\bullet}="5";
(80,3)*+{-}="-";
(80,5)*+{+}="+";
(100,0)*+{}="6";
(97,6)*+{}="6t";
(90,12)*+{}="6tt";
(97,-6)*+{}="6b";
(90,-12)*+{}="6bb";
 "5";"6" **\dir{-}; ?(.35)*\dir{<};
{\ar@{->} "5";"4"};
"5";"6t" **\dir{-}; ?(.35)*\dir{<};
"5";"6tt"  **\dir{-}; ?(.35)*\dir{<};
"5";"6b"  **\dir{-}; ?(.35)*\dir{<};
"5";"6bb"  **\dir{-}; ?(.35)*\dir{<};
\end{xy}
\end{eqnarray}
Signs $\pm$ over the vertices correspond to $\alpha^{\pm}$, while we neglect for simplicity the dependence on $z, \bar z$ variables. 
The incoming line in the vertex forms when during the normal ordering procedure the creation/annihilation operators from other vertices leave the commutator term \rf{lines} with $\alpha^{\pm}$ at a given vertex. The outgoing line forms when the creation/annihilation operator of a given vertex leaves a commutator term $\alpha^{\pm}$ from another vertex. 

On the other hand, $J^{\pm}$ also have terms of the form 
\begin{eqnarray}\label{term}
\kappa\p_u\alpha^{\pm}(z,\bar z), \quad \rho(z,\bar z)\alpha^{\pm}(z,\bar z).
\end{eqnarray}
According to the strategy formulated, we denote the first term from \rf{term} as a terminal vertex for graphs, because it contains only $\alpha^{\pm}$: 

\begin{eqnarray}
\begin{xy}
(60,0)*+{}="4";
(80,0)*+{\bullet}="5";
(80,3)*+{-}="-";
(80,5)*+{+}="+";
(100,0)*+{}="6";
(97,6)*+{}="6t";
(90,12)*+{}="6tt";
(97,-6)*+{}="6b";
(90,-12)*+{}="6bb";
 "5";"6" **\dir{-}; ?(.35)*\dir{<};
"5";"6t" **\dir{-}; ?(.35)*\dir{<};
"5";"6tt"  **\dir{-}; ?(.35)*\dir{<};
"5";"6b"  **\dir{-}; ?(.35)*\dir{<};
"5";"6bb"  **\dir{-}; ?(.35)*\dir{<};
\end{xy}
\end{eqnarray}

However the second term from \rf{term} is composite: it has contribution from  
$\alpha^{\pm}(z,\bar z)$ and $\rho(z,\bar z)$. We consider this term as a terminal vertex for the graphs coming from the normal ordering of $a,b$-operators, however, we should add one outgoing "wavy" line, corresponding to the normal ordering in the Fock space: 
 
 
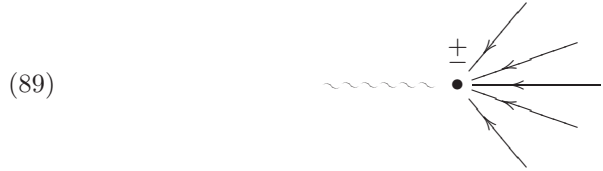
\begin{eqnarray}
\begin{xy}
(60,0)*+{}="4";
(80,0)*+{\bullet}="5";
(80,3)*+{-}="-";
(80,5)*+{+}="+";
(100,0)*+{}="6";
(97,6)*+{}="6t";
(90,12)*+{}="6tt";
(97,-6)*+{}="6b";
(90,-12)*+{}="6bb";
{\ar@{~} "5";"4"};
 "5";"6" **\dir{-}; ?(.35)*\dir{<};
"5";"6t" **\dir{-}; ?(.35)*\dir{<};
"5";"6tt"  **\dir{-}; ?(.35)*\dir{<};
"5";"6b"  **\dir{-}; ?(.35)*\dir{<};
"5";"6bb"  **\dir{-}; ?(.35)*\dir{<};
\end{xy}
\end{eqnarray}
 
Each wavy line produces the term from 2-point correlator \rf{2pointf} and may appear only once in the connected graph: it is so, because of the combinatorial formula, expressing the Fock space n-point correlator via 2-point correlators.

Finally, there are terms coming from $J^3$-generator, namely
\begin{eqnarray}
&&\label{ab2}\frac{1}{2}a(z,\bar z), \quad  \frac{1}{2}b(z,\bar z),\\  
&&\label{terminalj3}2\kappa\alpha^-(z,\bar z)\p_u \alpha^+(z,\bar z).
\end{eqnarray}
According to our conventions, we denote first two elements \rf{ab2} as initial vertices with one outgoing line to the right and to the left correspondingly:
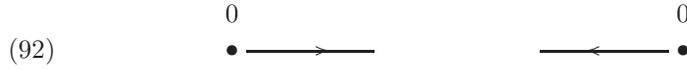
\begin{eqnarray}
\begin{xy}
(20,0)*+{\bullet}="2";
(20,5)*+{0}="0";
(40,0)*+{}="3";
"2";"3"  **\dir{-}; ?(.65)*\dir{>};
(60,0)*+{}="4";
(80,0)*+{\bullet}="5";
(80,5)*+{0}="0";
 "4";"5" **\dir{-}; ?(.35)*\dir{<};
\end{xy}
\end{eqnarray}
Zeros over the vertices represent the absence of $\alpha^{\pm}$ contribution from these terms.  
 Finally, \rf{terminalj3} will correspond to the terminal vertex, because it has neither creation or annihilation operators:  
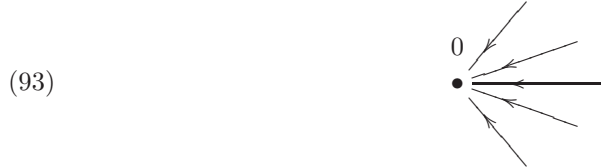
\begin{eqnarray}
\begin{xy}
(60,0)*+{}="4";
(80,0)*+{\bullet}="5";
(80,5)*+{0}="0";
(100,0)*+{}="6";
(97,6)*+{}="6t";
(90,12)*+{}="6tt";
(97,-6)*+{}="6b";
(90,-12)*+{}="6bb";
 "5";"6" **\dir{-}; ?(.35)*\dir{<};
"5";"6t" **\dir{-}; ?(.35)*\dir{<};
"5";"6tt"  **\dir{-}; ?(.35)*\dir{<};
"5";"6b"  **\dir{-}; ?(.35)*\dir{<};
"5";"6bb"  **\dir{-}; ?(.35)*\dir{<};
\end{xy}
\end{eqnarray}
However, in this exceptional case we will consider every line as a sum of two terms contributing $\delta(z,w)$, because there are two $\alpha$'s in this vertex corresponding to the term $\alpha^-(z,\bar z)\p_u\alpha^+(z,\bar z)$, so we consider contribution of the commutator with each of them. 
Moreover, formally on the circle (if such a limit exists), there can be only one incoming line into this vertex, 
 because of the commutation relations:
\begin{eqnarray}
[a(u), \alpha^+(v)\p_v \alpha^-(v)]=[b(u), \alpha^+(v)\p_v \alpha^-(v)]=
-\delta'(u-v).
\end{eqnarray}
Here are two sample diagrams 

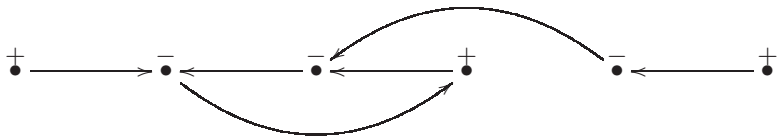
\begin{eqnarray}
\begin{xy}
(0,12)*+{+}="A+";
(20,12)*+{-}="B-";
(40,12)*+{-}="C-";
(60,12)*+{+}="D+";
(80,12)*+{-}="E-";
(0,10)*+{\bullet}="A";
(20,10)*+{\bullet}="B";
{\ar@{->} "A";"B"};
(40,10)*+{\bullet}="C";
(60,10)*+{\bullet}="D";
{\ar@{->} "C";"B"};
{\ar@/_2pc/ "B";"D"};
{\ar@{->} "D";"C"};
(80,10)*+{\bullet}="E";
{\ar@/_2pc/ "E";"C"};
(100,10)*+{\bullet}="F";
(100,12)*+{+}="F+";
{\ar@{->} "F";"E"}
\end{xy}\nonumber
\end{eqnarray}

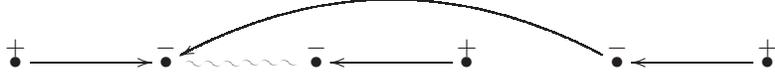
\begin{eqnarray}
\begin{xy}
(0,12)*+{+}="A+";
(20,12)*+{-}="B-";
(40,12)*+{-}="C-";
(60,12)*+{+}="D+";
(80,12)*+{-}="E-";
(0,10)*+{\bullet}="A";
(20,10)*+{\bullet}="B";
{\ar@{->} "A";"B"};
(40,10)*+{\bullet}="C";
(60,10)*+{\bullet}="D";
{\ar@{~} "C";"B"};
{\ar@{->} "D";"C"};
(80,10)*+{\bullet}="E";
(100,10)*+{\bullet}="F";
(100,12)*+{+}="F+";
{\ar@/_2pc/ "E";"B"}
{\ar@{->} "F";"E"}
\end{xy}\nonumber
\end{eqnarray}

which belong to the expansion of the 6-point correlator
\begin{eqnarray}
\langle J^+(z_1,\bar z_1)J^-(z_2,\bar z_2)J^-(z_3,\bar z_3)J^+(z_4,\bar z_4)J^-(z_5,\bar z_5)J^+(z_6,\bar z_6)\rangle.
\end{eqnarray}

Notice that every connected graph contributing to the correlator 
\begin{eqnarray}\label{corgen}
\langle \phi_1(z_1,\bar z_1)....\phi_n(z_n,\bar z_n)\rangle_p,
\end{eqnarray}
where $\phi_i=J^3, J^{\pm}$, has at most one loop. 
Actually, each vertex has possibly many incoming lines, but at most one outgoing line. Let us consider the graph with one loop. It means that all outgoing lines of the vertices in this loop are included in the loop and all other lines are incoming.  Therefore none of these vertices can participate in a different loop, since all outgoing lines are included in the first loop. Let us assume that there is another loop with different set of vertices in the same graph. However, since for all vertices participating all outgoing lines are "circulating" inside the loop there is no possibility to make a connected graph out of two loops. 

However, we see that each loop diagram, though  producing a well-defined expression for currents inside the circle, 
in the limit $|z_i|\to 1$ produce a divergence of the type $\delta(0)$, since every line in the loop produces delta-function. 
It is easy to see, say in the case of the correlator 
\begin{eqnarray}
\langle J^+(z_1,\bar z_1)J^-(z_2,\bar z_2)\rangle_p.
\end{eqnarray}
When the creation operator from $J^-$ contributes the commutator term \rf{lines} with the exponent in $J^+$ and at the same time, the annihilation operator from 
$J^-$ contributes the commutator term \rf{lines} with the exponent in $J^+$, the following diagram is produced: 
\begin{eqnarray}
\xymatrix{
\bullet \ar@/_1pc/[r] &
\bullet   \ar@/_1pc/[l]}
\end{eqnarray}
\\
\noindent and leads to the divergent term $\delta(u_1-u_2)\cdot\delta(u_1-u_2)$ when considered on the circle. 

The simplest way to regularize correlators to preserve commutation relations is to throw away all the loop diagrams.  This allows to make sense of the correlators on a circle and leads to the following Theorem \cite{fzsl2}. \\

\noindent {\bf Theorem 6.2.a.} {\it Let us consider the renormalized correlators
 \begin{eqnarray}\label{rencor}
\langle \phi_1(z_1,\bar z_1)....\phi_n(z_n,\bar z_n)\rangle_p^R,
\end{eqnarray}
where $\phi_i=J^3, J^{\pm}$, which contain only the tree graph contributions to \rf{corgen}.
Then the limit ${r_k\to 1}$ (where $z_k=r_ke^{iu_k}$) of \rf{rencor} exists as a distribution. Moreover, the commutation relation between $J^3, J^{\pm}$ under the renormalized correlator reproduces $\widehat{sl(2,\mathbb{R})}$. 
Therefore, the correlators \rf{rencor} considered on a circle define a module for $\widehat{sl(2,\mathbb{R})}$ with the Hermitian bilinear form.}\\

\noindent {\bf Proof.} 
To prove this theorem, it is enough to look at relevant contibutions for the commutators. It appears that they all come from initial/terminal vertices and their closest neighbors. For example, let us consider the diagram of this sort: 
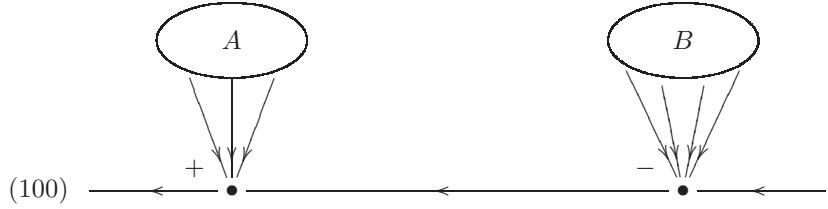
\begin{eqnarray}\label{diag1}
\begin{xy}
(0,0)*+{}="a";
(20,0)*+{\bullet}="b";
(80,0)*+{\bullet}="c";
(100,0)*+{}="d";
(15,3)*+{+}="+";
(75,3)*+{-}="-";
"a";"b"  **\dir{-}; ?(.45)*\dir{<};
"b";"c"  **\dir{-}; ?(.45)*\dir{<};
"c";"d"  **\dir{-}; ?(.45)*\dir{<};
%%%%%%%%%%%%%%%%%%%%%%%%%%%%%%%%
%(14,17)*+{}="b1";
%(20,15)*+{}="b2";
%(26,17)*+{}="b3";
%{\ar@{-} "b1";"b"};
%{\ar@{-} "b2";"b"};
%{\ar@{-} "b3";"b"};
(14,16)*+{}="b1";
(20,16)*+{}="b2";
(26,16)*+{}="b3";
"b1";"b"  **\dir{-}; ?(.75)*\dir{>};
"b2";"b"  **\dir{-}; ?(.75)*\dir{>};
"b3";"b"  **\dir{-}; ?(.75)*\dir{>};
%%%%%%%%%%%%%%%%%%%%%%%%%%%%%%%%%%%%%
(20,20)*+{A}="A";
(20,20)*\xycircle(10,5){};
(70,20)="e1";
(80,20)*+{B}="B";
(90,20)="e3";
(80,20)*\xycircle(10,5){};
%%%%%%%%%%%%%%%%%%%%%%%%%%%%%%%%%
%(72,17.5)*+{}="c1";
%(77,16)*+{}="c2";
%(83,16)*+{}="c3";
%(88,17.5)*+{}="c4";
%{\ar@{-} "c1";"c"};
%{\ar@{-} "c2";"c"};
%{\ar@{-} "c3";"c"};
%{\ar@{-} "c4";"c"};
(72,17)*+{}="c1";
(77,15)*+{}="c2";
(83,15)*+{}="c3";
(88,17)*+{}="c4";
"c1";"c"  **\dir{-}; ?(.75)*\dir{>};
"c2";"c"  **\dir{-}; ?(.72)*\dir{>};
"c3";"c"  **\dir{-}; ?(.72)*\dir{>};
"c4";"c"  **\dir{-}; ?(.75)*\dir{>};
\end{xy}
\end{eqnarray} 
\\
corresponding to the correlator 
\begin{eqnarray}\label{diag2}
\langle\dots J^+(z_1,\bar z_1)J^-(z_2,\bar z_2)\dots\rangle_p.
\end{eqnarray} 
There is a diagram with equal contribution 
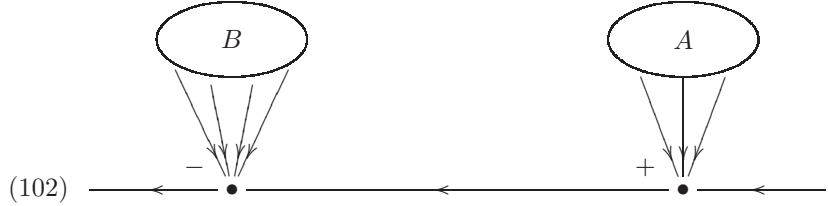
\begin{eqnarray}\label{diag2'}
\begin{xy}
(0,0)*+{}="a";
(20,0)*+{\bullet}="b";
(80,0)*+{\bullet}="c";
(100,0)*+{}="d";
(15,3)*+{-}="-";
(75,3)*+{+}="+";
"a";"b"  **\dir{-}; ?(.45)*\dir{<};
"b";"c"  **\dir{-}; ?(.45)*\dir{<};
"c";"d"  **\dir{-}; ?(.45)*\dir{<};
%%%%%%%%%%%%%%%%%%%%%%%%%%%%%%%%%%%%
%(12,17.5)*+{}="b1";
%(17,16)*+{}="b2";
%(23,16)*+{}="b3";
%(28,17.5)*+{}="b4";
%{\ar@{-} "b1";"b"};
%{\ar@{-} "b2";"b"};
%{\ar@{-} "b3";"b"};
%{\ar@{-} "b4";"b"};
(12,17)*+{}="b1";
(17,15)*+{}="b2";
(23,15)*+{}="b3";
(28,17)*+{}="b4";
"b1";"b"  **\dir{-}; ?(.75)*\dir{>};
"b2";"b"  **\dir{-}; ?(.72)*\dir{>};
"b3";"b"  **\dir{-}; ?(.72)*\dir{>};
"b4";"b"  **\dir{-}; ?(.75)*\dir{>};
%%%%%%%%%%%%%%%%%%%%%%%%%%%%%%%%%%%%
(20,20)*+{B}="B";
(20,20)*\xycircle(10,5){};
(70,20)="e1";
(80,20)*+{A}="A";
(90,20)="e3";
(80,20)*\xycircle(10,5){};
%%%%%%%%%%%%%%%%%%%%%%%%%%%%%%
%(74,17)*+{}="c1";
%(80,15)*+{}="c2";
%(86,17)*+{}="c3";
%{\ar@{-} "c1";"c"};
%{\ar@{-} "c2";"c"};
%{\ar@{-} "c3";"c"};
(74,16)*+{}="c1";
(80,16)*+{}="c2";
(86,16)*+{}="c3";
"c1";"c"  **\dir{-}; ?(.75)*\dir{>};
"c2";"c"  **\dir{-}; ?(.75)*\dir{>};
"c3";"c"  **\dir{-}; ?(.75)*\dir{>}
\end{xy}
\end{eqnarray}
\\
from the correlator
\begin{eqnarray}
\langle\dots J^-(z_2,\bar z_2)J^+(z_1,\bar z_1)\dots\rangle_p.
\end{eqnarray}
Therefore, they cancel each other when we take a commutator of $J^+(u_1), J^-(u_2)$, since the middle line produces $\delta(u_1-u_2)$. One can understand it in the following way: the commutator 
$[J^+(u_1), J^-(u_2)]$ involves $J^3$-term and 
$\delta'$-term, 
$J^3$-term contains only initial and terminal vertices, therefore the diagrams, which contribute to this commutator should have also $J^+$ or $J^-$ as initial or terminal vertices. As for the commutators 
of $J^3$ with $J^{\pm}$, the picture above indicates that only terminal/initial vertices and their closest neighbors contribute, because again, all terms in $J^3$ are depicted as initial/terminal vertices. Therefore, it is clear that one can cancel all loop contributions to the correlators and still the commutation relations of Theorem 6.1 will be valid. However, if all the loops are eliminated, one can consider the limit $|z_i|\to 1$, putting all currents on a circle, since all divergent graphs are gone. Therefore, we have the well-defined correlators of the generators of $\widehat{sl(2,\mathbb{R})}$ Lie algebra. One can see  
that these correlators define a Hermitian bilinear form, because after the conjugation one can obtain one-to-one correspondence between tree graphs in conjugated correlators.
\hfill$\blacksquare$\bigskip
\\

It appears that the theorem above can be generalized: instead of eliminating loops completely, one can renormalize them, i.e. eliminate the divergence of the type $\delta(0)$. Namely, one can associate a real number $\mu_k$ with every loop which is the contribution of $k$ $J^{\pm}$-currents at the points $z_1,..., z_k$. This number $\mu_k$ enters the loop in the following way. 
Suppose currents $J^{\pm}(z_i)$ appear in the correlator in the given order from left to right. Our loop contains the following term
\begin{eqnarray}\label{delta1}
\delta (z_1,z_2)\cdot\delta(z_2,z_3)\dots\delta(z_{k-1}, z_k)\cdot \delta(z_{k}, z_1),
\end{eqnarray}
where $\delta (z,w)=\sum_{n\ge 0}(z\bar w)^n+\sum_{n< 0}(w\bar z)^n$. In the limit $|z_k|\to 1$ we will substitute it with  
the expression 
\begin{eqnarray}\label{delta2}
\mu_k\cdot \delta(u_1-u_2)\cdot\delta(u_2-u_3)\dots\delta(u_{k-1}-u_k),
\end{eqnarray}
where we remind that $z_k=r_ke^{iu_k}$. 
This procedure again does affect neither commutation relations inside the correlator nor Hermicity condition for the resulting bilinear form, because loops do not participate in commutation relations. Therefore,  we have a new version of the theorem above \cite{fzsl2}.
\\

\noindent {\bf Theorem 6.2.b.} {\it Let us consider the renormalized correlators
 \begin{eqnarray}\label{rencormu}
\langle \phi_1(z_1,\bar z_1)....\phi_n(z_n,\bar z_n)\rangle_p^{R, \{\mu_n\}},
\end{eqnarray}
where $\phi_i=J^3, J^{\pm}$, such that the loop contributions to \rf{corgen} are replaced by their renormalized analogues with the family of arbitrary real parameters $\{\mu_n\}$. 
Then the limit ${r_k\to 1}$ (where $z_k=r_ke^{iu_k}$) of \rf{rencormu} exists and the commutation relations between $J^3, J^{\pm}$ under the renormalized correlator reproduce Lie algebra $\widehat{sl(2,\mathbb{R})}$, and correlators \rf{rencormu} considered on a circle define a module for $\widehat{sl(2,\mathbb{R})}$ with the Hermitian bilinear form.}\\

\noindent {\bf Proof.} The proof goes along the same lines as in the part $a$ of the theorem. When we consider the loop contribution to the commutator, there will be two loop diagrams, each emerging from one of the terms precisely as in the pictures \rf{diag1}, \rf{diag2} which again will cancel each other, because of the delta function line between two vertices and because there always is one incoming line and outgoing line for each vertex associated with the currents, participating in the commutator in these two terms. Therefore, the commutation relations are unaffected in the presence of loops. It is easy to see that the same applies to the Hermicity condition.
\hfill$\blacksquare$\bigskip
\\

We mention that in this article the construction of $\widehat{sl}(2,\mathbb{R})$-modules based on unitary modules of $\mathcal{K}$-algebra is not considered. The renormalization procedure in that case is very similar, the difference is that there are several tree diagrams in teh expansion of correlators, which also require renormalization procedure. For more details on the treatment of unitary representations of unitary modules of $\mathcal{K}$-algebra one can consult \cite{fzsl2}.

Finally, it is important to remark that the realization of $\widehat{sl(2,\mathbb{R})}$ from Section 6.1 can be generalized in the higher rank case via multiple $\mathcal{K}$-algebras leading to the formulas similar to \rf{sl2R1}, \rf{sl2R2} (in the case of quantum group there is a close construction using algebras of quantum planes \cite{ivan2}, \cite{ivan3}). It will be a pure technical problem to modify Theorem 6.2.a and 6.2.b in the case of a  higher rank.

\section*{Appendix: Classical Wiener Measure and Gaussian Integration on Hilbert spaces}

\noindent{\bf A.1. A brief review of the abstract approach to Wiener measure.} In this section we collect all the necessary facts about the Wiener measure. For a more detailed exposition of this subject one can consult \cite{quo}, \cite{hida}, \cite{gelfand}, \cite{lee}, \cite{frenkel}. 

The abstract Wiener measure is the Gaussian measure on a Banach space with certain properties. 
The construction is as follows. We start from the Hilbert space $\mathcal{H}$ and consider a Gaussian measure associated with the unital operator and the norm $||\cdot||$ on the Hilbert space, such that heuristically measure can be represented as follows:
\begin{eqnarray} 
d\tilde{w}^t\sim e^{-\frac{||x||^2}{2t}}[dx].
\end{eqnarray}
However, this Gaussian measure is not $\sigma$-additive. In order to make it $\sigma$-additive, one  has to 
consider a weaker norm $|\cdot|$ with certain conditions on it with respect to measure $d\tilde{w}^t$. Then one can consider a completion $\mathcal{H}$ with respect to the norm $| \cdot |$. This will give a Banach space $\mathcal{B}$. Now the Gaussian measure $d\tilde w^t$ can be extended to the Banach space, where it becomes 
$\sigma$-additive. 

The important example, which we will consider in the next subsection, is constructed as follows. Let us take the Hilbert space $C_0'[0,2\pi]$ of real-valued absolutely continuous functions, such that $x(0)=0$ for any $x\in C_0'[0,2\pi]$ and 
$\int_0^{2\pi}(x'(u))^2du<\infty$. The inner product is given by: $\langle x_1,x_2\rangle=\int_0^{2\pi}x_1'(u)x_2'(u)du$. 
%One can notice, that the derivative 
%here creates an isometry map from $C_0'[0,2\pi]$ to the Hilbert space of square integrable real-valued functions on $[0,2\pi]$, i.e. $L^2([0,2\pi];\mathbb{R})$. 
Then one can consider the weaker norm on this space: $|x|={\rm sup}_{u\in[0,2\pi]} x(u)$. The completion of $C_0'[0,2\pi]$ with respect to $|\cdot|$ is the Banach space $C_0[0,2\pi]$ of real-valued continuous functions such that $x(0)=0$ for any $f\in C_0[0,2\pi]$. 

It appears that the norm $|\cdot|$ satisfies all necessary properties and therefore there exists a $\sigma$-additive Gaussian measure $dw^t$ on $C_0[0,2\pi]$. This measure is called the classical Wiener measure. In the next subsection we will give its direct construction, using another approach.\\

\noindent {\bf A.2. The construction of the classical Wiener measure and its basic properties.} Consider the space of continuous functions $C_0[0,2\pi]$. Let $C_{0,X}[0,2\pi]$ denote the closed subspace of  $C_0[0,2\pi]$, consisting 
of real-valued continuous functions such that $x(2\pi)=X$ for some $X\in \mathbb{R}$. 

The following subsets of $C_0[0,2\pi]$ are called cylinder sets:
\begin{eqnarray} 
\{x\in C_0[0,2\pi]: x(\tau_1)\in A_1, \dots , x(\tau_n) \in A_n, 0<\tau_1<\dots\tau_n\le 2\pi\},
\end{eqnarray}
where $A_1,A_2,\dots A_n$ are Borel subsets of $\mathbb{R}$.\\
The Wiener measure with variance $t\ge 0$ is defined on the cylindrical sets of $C_0[0,2\pi]$ as follows:
\begin{eqnarray}
&&w^t( x(\tau_1)\in A_1, \dots , x(\tau_n) \in A_n)=\nonumber\\
&&\int_{A_1}\dots \int_{A_n} u_t(\Delta x_1, \Delta \tau_1)\dots 
u_t(\Delta x_n, \Delta \tau_n)dx_1\dots dx_n.
\end{eqnarray}
The conditional Wiener measure of variance $t>0$ is defined on the cylinder sets of $C_{0,X}[0,2\pi]$ as follows:
\begin{eqnarray}
&&w^t_X( x(\tau_1)\in A_1, \dots , x(\tau_n) \in A_n)=\nonumber\\
&&\int_{A_1}\dots \int_{A_{n-1}} u_t(\Delta x_1, \Delta \tau_1)\dots 
u_t(\Delta x_n, \Delta \tau_n)dx_1\dots dx_{n-1},
\end{eqnarray}
where $u_t(x,s)=\frac{1}{2\pi\sqrt {ts}}e^{-\frac{x^2}{4\pi st}}$, $dx$ stands for the Lebesgue measure on $\mathbb{R}$, 
$\Delta x_k=x_k-x_{k-1}$, $\Delta \tau_k=\tau_k-\tau_{k-1}$, $x_0=0$, $\tau_0=0$. In the case ii) $x_n=X$, $\tau_n=2\pi$.

It appears that the resulting measures are $\sigma$-additive and moreover  $w^t(C_0[0,2\pi])=1$, $w^t_X(C_{0,X}[0,2\pi])=u_t(X,2\pi)$. 

One can show that the description of the Wiener measure that we gave in subsection A.1. and the one introduced in this subsection agree (see Section I.5. of \cite{quo}). 
Now we want to relate $L^2$ spaces with respect to the Wiener measure and the conditional Wiener measure. In order to do that, notice the following property. If $f$ is an integrable function on $C_{0,X}([0,2\pi])$ for almost all $X$, then by the definition of the Wiener measure one has:
\begin{eqnarray}
\int_{C_{0}[0,2\pi]} f(x)dw^t(x)=\int_{\mathbb{R}}\Big(\int_{C_{0,X}[0,2\pi]}f(x)dw_X^t(x)\Big)dX.
\end{eqnarray}
Therefore, the Hilbert space of square-integrable functions with respect to $dw^t$ decomposes as a direct integral:
\begin{eqnarray}
L^2(C_{0}[0,2\pi], dw^t)=\int_{\mathbb{R}}^{\oplus}L^2(C_{0,X}([0,2\pi]), dw^t_X)dX.
\end{eqnarray}

The Wiener measure and the conditional Wiener measure have translation property. Let us describe this property in detail.

Let $f$ be an integrable function on  $C_{0}([0,2\pi])$ and $y\in C'_0[0,2\pi]$, then 
\begin{eqnarray}
&&\int_{C_{0}[0,2\pi]} f(x)dw^t(x)=\\
&&\int_{C_{0}[0,2\pi]} f(x+y)e^{-\frac{1}{t}\int^{2\pi}_0y'(u)dx(u)-\frac{1}{2t}\int^{2\pi}_0y'(u)y'(u)du}dw^t(x),\nonumber
\end{eqnarray}
where $\int^{2\pi}_0y'(u)dx(u)$ is the Stieltjes integral.\\

If $f$ is an integrable function in $C_{0,X}([0,2\pi])$ and $y$ as above, such that $y(2\pi)=Y$ then,
\begin{eqnarray}
&&\int_{C_{0,X+Y}[0,2\pi]} f(x)dw_{X+Y}^t(x)=\nonumber\\
&&\int_{C_{0,X}[0,2\pi]} f(x+y)e^{-\frac{1}{t}\int^{2\pi}_0y'(u)dx(u)-\frac{1}{2t}\int^{2\pi}_0y'(u)y'(u)du}dw_{X}^t(x).
\end{eqnarray}

The space of Wiener measure translations, i.e. $C'_0[0,2\pi]$ is usually called the Cameron-Martin space.

 One can define a unitary operator on $L^2$ space for any abstract Wiener measure, which is similar to the Fourier transform. We will need it only in the case of measure $dw^t_0$ on $C_{0,0}[2\pi]$.
The formula is as follows:
\begin{eqnarray}
Ff(y)=\int_{C_{0,0}[0,2\pi]}f(x+iy)dw^{2t}_0(x).
\end{eqnarray} 
Here $f\in L^2(C_{0,0}[0,2\pi], dw^t_0)$. The inverse transformation is given by the changing sign of $y$ in the formula 
above.\\

\noindent{\bf A.3. Gaussian integration on Hilbert spaces.} 
In this subsection we recall a few basic facts and formulas, for more information, see e.g. \cite{daprata}, \cite{quo}. 
Suppose we have a real separable Hilbert space $H$ with the orthonormal basis $\{e^i\}$, $i\in \mathbb{N}$ and the pairing $\langle\cdot , \cdot\rangle$. Every element $x$ of this Hilbert space can be expressed as $x=\sum_i x_ie^i$. Let us introduce positive real numbers $\lambda_i$, $i\in \mathbb{N}$, so that 
 $\sum_i\lambda_i<\infty$. Then one can say that the numbers $\lambda_i$ define a diagonal trace class operator $A$ on our Hilbert space. 
Than it appears possible to define the sigma-additive measure $d\mu_A$ on $H$ 
and heuristically express it as follows:
\begin{eqnarray}
d\mu_A(x)=(\sqrt{\det {2 \pi A}})^{-1}\cdot e^{-\frac{1}{2}\langle x,A^{-1}x \rangle}[dx],
\end{eqnarray}
 which can be thought as the infinite product of 1-dimensional Gaussian measures for each $i$: $d\mu_i=\sqrt{2\pi \lambda_i}^{-1}e^{-\lambda_i^{-1}x_i^2}$.

Since it is a sigma-additive measure, one can define the space of square-integrable functions $L^2(H, d\mu_A)$ with respect to it. One of the basic formulas
is the translational shift in the measure. Namely, if $b\in Im A$, then 
\begin{eqnarray}
\int f(x) d\mu_A(x)=\int f(x+b) e^{{-\frac{1}{2}\langle b,A^{-1}b \rangle}-\langle x,A^{-1}b \rangle}d\mu_A(x).
\end{eqnarray}
One can consider also an infinitesimal version of this formula. Making $b$ infinitesimal and parallel to $e_i$, we obtain 
that 
\begin{eqnarray}
\int D_i f(x)d\mu_A(x)=0, \quad D_i=\p_{x_i}-\lambda_i^{-1}x_i,
\end{eqnarray}
if $\p_{x_i} f(x)\in L^2(H, d\mu_A)$.  
It should be noted that the following monomials
\begin{eqnarray}\label{polexp}
(\prod^k_{i=1}\langle \alpha_i, x\rangle) e^{\langle \beta, x\rangle},
\end{eqnarray}
where $\alpha_i, \beta$ are the elements of the complexified Hilbert space, are always integrable with respect to $d\mu_A$, moreover, they belong  to $L^2(H, d\mu_A)$. There is an explicit formula for the integral of the function \rf{polexp}, which can be derived from the simple result:
\begin{eqnarray}\label{gauss}
\int  e^{\langle \beta, x\rangle}d\mu_A(x)=e^{\frac{1}{2}\langle \beta, A\beta\rangle}.
\end{eqnarray}

\end{document}